\documentclass{article}

\input{diagrams}
\diagramstyle[repositionpullbacks,midshaft,PostScript=dvips,nohug,scriptlabels,size=2em] 
\usepackage{graphics}
\usepackage{color}
\usepackage[mathcal,mathscr]{euscript}

\usepackage{amssymb}
\newcommand{\emptybk}{\:\:}

\newcommand{\dashbk}{-}
\newcommand{\diagspace}{\mbox{\hspace{2em}}}
\newcommand{\bref}[1]{(\ref{#1})}
\newcommand{\ucontents}[2]{\addcontentsline{toc}{#1}{\numberline{}{#2}}}
\newcommand{\mb}[1]{\mathbf{#1}}

\newcommand{\mr}[1]{\mathrm{#1}}

\newcommand{\fcat}[1]{\mb{#1}}
\newcommand{\ovln}[1]{\overline{#1}}
\newcommand{\twid}[1]{\widetilde{#1}}
\newcommand{\slsh}{/\linebreak[0]}
\newcommand{\dt}{.\linebreak[0]}
\newcommand{\such}{\:|\:}
\newcommand{\without}{\setminus}
\newcommand{\blob}{{\raisebox{.3ex}{\ensuremath{\scriptscriptstyle{\bullet}}}}}
\newcommand{\epsln}{\varepsilon}
\newcommand{\implies}{\,\Rightarrow\,}
\newcommand{\op}{\mr{op}}
\newcommand{\id}{\mr{id}}

\newcommand{\One}{\fcat{1}}
\newcommand{\Set}{\fcat{Set}}
\newcommand{\Top}{\fcat{Top}}
\newcommand{\ftrcat}[2]{[#1,#2]}
\newcommand{\goesto}{\,\longmapsto\,}
\newcommand{\goiso}{\goby{\diso}}

\newcommand{\oppair}[2]{\pile{\rTo^{\scriptstyle #1}\\ 
\lTo_{\scriptstyle #2}}} 
\newcommand{\parpair}[2]{\pile{\rTo^{\scriptstyle #1}\\ 
\rTo_{\scriptstyle #2}}}

\newarrow{Get}....>
\newarrow{Goesto}|---{->}
\newarrow{Incl}C--->
\newarrow{Mod}--+->
\newcommand{\diso}{\sim}

\newcommand{\upr}[1]{[#1]}	
\newcommand{\reals}{\mathbb{R}}
\newcommand{\latin}[1]{\textit{#1}}
			
\newenvironment{scriptarray}%
{\renewcommand{\arraystretch}{0.8}\begin{array}{l}}%
{\end{array}\renewcommand{\arraystretch}{1}}
\newcommand{\demph}[1]{\textbf{\textup{#1}}}
\newcommand{\done}{\hfill\ensuremath{\Box}}
\newenvironment{prooflike}[1]{\begin{trivlist}\item\textbf{#1}\ }
{\end{trivlist}}
\newenvironment{proof}{\begin{prooflike}{Proof}}{\end{prooflike}}
\newcommand{\scat}[1]{\mathbb{#1}}
\newcommand{\go}{\rTo\linebreak[0]}
\newcommand{\goby}[1]{\rTo^{#1}\linebreak[0]}
\newcommand{\iso}{\cong}
\newcommand{\nat}{\mathbb{N}}	
\newcommand{\url}[1]{#1}

\newcommand{\elt}[1]{\scat{E}\left(#1\right)}
\newcommand{\gomod}{\rMod\linebreak[0]}
\newcommand{\gobymod}[1]{\rMod^{#1}\linebreak[0]}

\newcommand{\cstyle}[1]{\textbf{\upshape{#1}}}
\newcommand{\So}{\cstyle{S}}

\newarrow{Qt}---->
\newarrow{Modget}..+.>

\newcommand{\frc}[2]{\frac{#1}{#2}}
\newcommand{\half}{\frc{1}{2}}
\newcommand{\dfrc}[2]{\textstyle{\frac{#1}{#2}}}
\newcommand{\dhalf}{\dfrc{1}{2}}

\newcommand{\ob}{\mr{ob}}

\newcommand{\diam}{\mr{diam}}

\newcommand{\mod}{\mr{mod}\ } 
\newcommand{\Dface}{\Delta_{\mr{face}}}
\newcommand{\propersub}{\subset}
\newcommand{\nepower}{\mathcal{P}_{\neq\emptyset}}
\newcommand{\vtr}[1]{#1}
\newcommand{\oba}{\mr{ob}\,}

\newcommand{\restr}[2]{#1|_{#2}}
\newcommand{\nem}[1]{#1^+}
\newcommand{\of}{\,\raisebox{0.08ex}{\ensuremath{\scriptstyle\circ}}\,}

\newcommand{\sub}{\subseteq}
\newcommand{\catI}{\mathcal{I}}

\newenvironment{diagdiag}
	{\begin{diagram}[size=1.5em]}
	{\end{diagram}}
\newcommand{\cell}[4]{\put(#1,#2){\makebox(0,0)[#3]{\ensuremath{#4}}}}
\newcommand{\zmark}{\scriptstyle{\bullet}}
\definecolor{grey}{gray}{0.5}
\newcommand{\avdots}{\rotatebox{90}{$\cdots$}}

\newtheorem{thm}{Theorem}[section]
\newtheorem{propn}[thm]{Proposition}
\newtheorem{lemma}[thm]{Lemma}
\newtheorem{cor}[thm]{Corollary}
\newtheorem{lotsofremarks}[thm]{Remarks}

\newtheorem{conj}[thm]{Conjecture}

\newtheorem{predefn}[thm]{Definition}

\newtheorem{preexample}[thm]{Example}
\newenvironment{example}{\begin{preexample}\upshape}{\end{preexample}}
\newenvironment{example*}[1]{\begin{preexample}[#1]\upshape}{\end{preexample}}
\newtheorem{prewarning}[thm]{Warning}

\newcommand{\xrefstyle}[1]{#1}

\newcommand{\secselfsimilaritysystems}{\xrefstyle{1}}
\newcommand{\secunivsoln}{\xrefstyle{4}}

\newcommand{\thmtopologicalFreyd}{\xrefstyle{1.7}}
\newcommand{\defnMcoalgebrainTop}{\xrefstyle{1.9}}
\newcommand{\corcategorywithpullbacksandequalizers}{\xrefstyle{2.9}}
\newcommand{\egnondegencofork}{\xrefstyle{2.12}}
\newcommand{\lemmaequalityinMtensorX}{\xrefstyle{3.2}}
\newcommand{\lemmaKoenig}{\xrefstyle{5.1}}
\newcommand{\lemmaequalityinIa}{\xrefstyle{6.2}}
\newcommand{\lemmarelationsdetermineequality}{\xrefstyle{6.3}}
\newcommand{\lemmasequentialcompactness}{\xrefstyle{6.5}}
\newcommand{\propnIaasaquotient}{\xrefstyle{6.7}}
\newcommand{\lemmacoprojectionsclosed}{\xrefstyle{6.14}}
\newcommand{\thmexistenceofuniversalsolution}{\xrefstyle{A.1}}

\newcommand{\eqinformalsss}{\xrefstyle{6}}
\newcommand{\eqFreydSSS}{\xrefstyle{9}}

\begin{document}

\sloppy

\title{A general theory of self-similarity II:\\ 
recognition}
\author{Tom Leinster}
\date{\normalsize 
Department of Mathematics, University of Glasgow\\
www\dt maths\dt gla\dt ac\dt uk\slsh $\sim$tl\\
tl@$\,\!$maths\dt gla\dt ac\dt uk}

\maketitle
\thispagestyle{empty}

\begin{center}
\textbf{Abstract}
\end{center}
This paper concerns the self-similarity of topological spaces, in the sense
defined in~\cite{SS1}.  I show how to recognize self-similar spaces, or
more precisely, universal solutions of self-similarity systems.  Examples
include the standard simplices (self-similar by barycentric subdivision)
and solutions of iterated function systems.  Perhaps surprisingly, every
compact metrizable space is self-similar in at least one way.  From this
follow the classical results on the role of the Cantor set among compact
metrizable spaces.

\vfill

\tableofcontents

\section*{Introduction}
\ucontents{section}{Introduction}

This paper's predecessor~\cite{SS1} introduced the notions of
`self-similarity system' (to be thought of as a system of simultaneous
equations in which each member of a family of objects is expressed as a
gluing-together of other members) and `universal solution' of such a
system.  A self-similarity system possesses a universal solution if and
only if an explicit condition \So\ holds; if so then the universal solution
is unique and can be constructed explicitly.

Few examples have been given so far.  In principle one can take any
self-similarity system $(\scat{A}, M)$ satisfying \So\ and find the
universal solution by going through the explicit construction.  In practice
this is cumbersome and it is much quicker to apply one of the Recognition
Theorems proved in \S\ref{sec:recognitiontheorems} below, as follows.

Recall that if $(J, \gamma)$ is a universal solution of $(\scat{A}, M)$
then $J \iso M \otimes J$, that is, each space $Ja$ must be homeomorphic to
the prescribed gluing-together of $(Jb)$s.  This necessary condition is not
sufficient (see \cite{SS1}, after Definition~\defnMcoalgebrainTop), but it
turns out to become sufficient if, roughly speaking, each space $Ja$ is
nonempty, compact, and made up of \emph{smaller} copies of $(Jb)$s.  Put
another way, a solution $(J, \gamma)$ of the system is universal if it
consists of nonempty compact metric spaces and contractions between them.
This is the Crude Recognition Theorem, essentially; it gives sufficient
conditions for an $M$-coalgebra to be a universal solution.  A more refined
result is the Precise Recognition Theorem, which gives necessary and
sufficient conditions.

To find the universal solution of a given self-similarity system, one can
therefore make a guess and check that it satisfies the hypotheses of a
Recognition Theorem.  I do this in~\S\ref{sec:examples} for various
specific self-similarity systems, formulating, for instance, the
self-similarity of the cubes $[0, 1]^n$, the simplices $\Delta^n$ (by both
barycentric and edgewise subdivision), the Cantor set, and the solutions of
certain iterated function systems.

A further challenge is to recognize which topological spaces arise as
universal solutions to some self-similarity system.  The Recognition
Theorems show that all such spaces are compact and metrizable, so the
question is: which compact metrizable spaces are self-similar?  Perhaps
surprisingly, the answer turns out to be: all of them~(\S\ref{sec:sass}).

A measure of the non-triviality of this result is that some of the
classical results on the Cantor set can be derived as corollaries: for
instance, that every nonempty compact metric space is a quotient of the
Cantor set.  This is the subject of~\S\ref{sec:thecantorset}.

Terminology and notation from~\cite{SS1} are used without further
mention.  An overview of~\cite{SS1} and the present paper is
available:~\cite{GSSO}. 

I gratefully acknowledge a Nuffield Foundation award NUF-NAL 04.

\section{Recognition theorems}
\label{sec:recognitiontheorems}

Fix a self-similarity system $(\scat{A}, M)$.

If $(\scat{A}, M)$ has a universal solution then by
Theorem~\thmexistenceofuniversalsolution\ of~\cite{SS1}, condition \So\
holds and the universal solution is the $(I, \iota)$ constructed in
\S\secunivsoln\ of \cite{SS1}.  We begin by establishing various properties
of $(I, \iota)$, working throughout in $\Top$ rather than $\Set$.

First, Lambek's Lemma guarantees that $(I, \iota)$ is a \demph{fixed point}
of $M$, that is, an $M$-coalgebra whose structure map is an isomorphism.
(Note that by definition, fixed points are nondegenerate.)  A fixed point
$(J, \gamma)$ is a coalgebra, but can also be regarded as an algebra $(J,
\psi)$ where $\psi = \gamma^{-1}$.  By definition, an \demph{$M$-algebra}
(in $\Top$) is a nondegenerate functor $J: \scat{A} \go \Top$ together with
a map $\psi: M \otimes J \go J$; but
\[
(M \otimes J) a 
=
\int^b M(b, a) \times Jb,
\]
so $\psi$ consists of a map $\psi_m: Jb \go Ja$ for each module element $m:
b \gomod a$, satisfying naturality axioms.  These say that if
\[
b' \goby{g} b \gobymod{m} a \goby{f} a'
\]
then $\psi_{fmg} = (Jf) \of \psi_m \of (Jg)$.  

For example, the fixed point $(I, \iota)$ has algebra structure $\phi =
\iota^{-1}$, where the components $\phi_m$ are as defined in~\cite{SS1}
(after Lemma~\lemmaequalityinIa).

\begin{lemma}
\label{lemma:fixedpointcomponents}
Let $(J, \psi)$ be a fixed point of $M$.  Then for each module element $b
\gobymod{m} a$, the map $Jb \goby{\psi_m} Ja$ is closed.
\end{lemma}

\begin{proof}
$\psi_m$ is the composite
\[
Jb
\goby{m \otimes \dashbk}
(M \otimes J) a
\rTo^{\psi_a}_\diso
Ja
\]
and $m \otimes \dashbk$ is closed by~\cite[\lemmacoprojectionsclosed]{SS1}
and nondegeneracy of $J$.
\done
\end{proof}

Second, the sets $Ia$ are not empty unless they are forced to be.  If, for
instance, $a \in \scat{A}$ and there is no module element of the form $b
\gobymod{m} a$, then for any coalgebra $(X, \xi)$ we have $\xi_a: Xa \go (M
\otimes X)a = \emptyset$, so $Xa$ is forced to be empty.  More generally,
$Xa$ is forced to be empty if there is no infinite chain $\cdots \gomod a_1
\gomod a$.  Say that a functor $X: \scat{A} \go \Set$ is \demph{occupied}
if $\catI a \neq \emptyset \implies Xa \neq \emptyset$ for all $a \in
\scat{A}$.  Trivially, $I$ is occupied.

Third, the spaces $Ia$ are metrizable, by
\begin{lemma}	\label{lemma:compactmetrizable}
A compact space is metrizable if and only if it is Hausdorff and
second countable.
\end{lemma}
\begin{proof}
See~\cite[IX.2.9]{Bou2} (where `compact' means compact Hausdorff).
\done
\end{proof}

One naturally asks how a metric can be defined.  There are many
possible metrics and apparently no canonical choice among them,
but the following result tells us all we need to know.  Recall
from~\cite{SS1} (before Lemma~\lemmarelationsdetermineequality)
that for each $a \in \scat{A}$ and $n\in\nat$ we have a closed
binary relation $R_n^a$ on $Ia$.

\begin{propn}[Metric on $Ia$]
\label{propn:metriconIa}
For each $a \in \scat{A}$ there exists a metric $d$ on $Ia$
compatible with the topology on $Ia$ and such that
\[
\forall \epsln > 0,
\ 
\exists n \in \nat:
\ 
(t, t') \in R_n^a
\ 
\implies 
\ 
d(t, t') \leq \epsln.
\]
\end{propn}
\begin{proof}
In fact, the stated property holds for any compatible metric $d$
on $Ia$ (and by Lemma~\ref{lemma:compactmetrizable}, there is at
least one such).  For if it does not hold, there is some $\epsln
> 0$ such that $d^{-1} [\epsln, \infty) \cap R_n^a \neq
\emptyset$ for all $n\in\nat$, where $d^{-1} [\epsln, \infty)$ is
the inverse image of $[\epsln, \infty)$ under $d: Ia \times Ia
\go \reals$.  Since $d$ is continuous, $d^{-1} [\epsln, \infty)$
is closed.  Each $R_n^a$ is also closed and $Ia$ is compact, so
$d^{-1} [\epsln, \infty) \cap \bigcap_{n\in\nat} R_n^a \neq
\emptyset$; but by Lemma~\lemmarelationsdetermineequality\ of
\cite{SS1}, $\bigcap_{n\in\nat} R_n^a = \Delta_{Ia}$, a
contradiction.  \done
\end{proof}

\begin{thm}[Precise Recognition Theorem]
\label{thm:preciserecognition}
The following are equivalent conditions on a fixed point $(J, \gamma =
\psi^{-1})$ of $M$ in $\Top$:
\begin{enumerate}
\item	\label{item:pr-univ}
$(J, \gamma)$ is a universal solution of $(\scat{A}, M)$ in $\Top$
\item	\label{item:pr-met}
$J$ is occupied, and for each $a \in \scat{A}$ the space $Ja$ is compact
and can be metrized in such a way that
\[
\renewcommand{\arraystretch}{1.2}
\begin{array}{l}
\forall \epsln > 0,
\ 
\exists n \in \nat:
\ 
\forall (a_n \gobymod{m_n} \cdots \gobymod{m_1} a_0) \in \catI_n a,
\\
\diam (\psi_{m_1} \cdots \psi_{m_n} (Ja_n)) \leq \epsln
\end{array}
\renewcommand{\arraystretch}{1}
\]
\item	\label{item:pr-notmet}
$J$ is occupied, and for each $a \in \scat{A}$, the space $Ja$ is compact
and 
\[
\forall (\cdots \gobymod{m_1} a_0) \in \catI a,
\ 
\left|
\bigcap_{n\in\nat}
\psi_{m_1} \cdots \psi_{m_n} (Ja_n)
\right|
\leq
1
\]
(where $| \emptybk |$ means cardinality).
\end{enumerate}
\end{thm}

\begin{proof}
\paragraph*{\bref{item:pr-univ}$\implies$\bref{item:pr-met}}
Assume~\bref{item:pr-univ}.  By Theorem~\thmexistenceofuniversalsolution\
of \cite{SS1}, condition~\So\ holds, so $(J, \gamma)$ is the universal
solution $(I, \iota)$ constructed in~\S\secunivsoln\ of~\cite{SS1}.
Certainly $I$ is occupied and each $Ia$ is compact.
Proposition~\ref{propn:metriconIa} then gives metrics with exactly the
property required.

\paragraph*{\bref{item:pr-met}$\implies$\bref{item:pr-notmet}}
Trivial.

\paragraph*{\bref{item:pr-notmet}$\implies$\bref{item:pr-univ}}
Assume~\bref{item:pr-notmet}.  First I claim that each $Ja$ is Hausdorff.
For each $n \in \nat$, put 
\[
\twid{R}_n^a 
=
\bigcup
\{
(\psi_{m_1} \cdots \psi_{m_n} (Ja_n))^2
\such
(a_n \gobymod{m_n} \cdots \gobymod{m_1} a_0)
\in 
\catI_n a
\}
\sub
Ja \times Ja.
\]
Then each $\twid{R}_n^a$ is closed by
Lemma~\ref{lemma:fixedpointcomponents} and finiteness of $\catI_n a$, so it
suffices to show that $\bigcap_{n\in\nat} \twid{R}_n^a = \Delta_{Ja}$.
Certainly $\bigcap_{n\in\nat} \twid{R}_n^a \supseteq \Delta_{Ja}$.
Conversely, let $(t, t') \in \bigcap_{n\in\nat} \twid{R}_n^a$.  For each
$n\in\nat$ we may choose
\[
(
a^n_n \gobymod{m^n_n} \cdots \gobymod{m^n_1} a^n_0
)
\in 
\catI_n a
\]
such that $t, t' \in \psi_{m^n_1} \cdots \psi_{m^n_n} (J a^n_n)$.  By
K\"onig's Lemma~\cite[\lemmaKoenig]{SS1} and finiteness of $\catI_n a$, we
may then choose
\[
(
\cdots \gobymod{m_2} a_1 \gobymod{m_1} a_0
)
\in
\catI a
\]
such that for all $r \in \nat$, there exists $n \geq r$ satisfying
\[
(a_r \gobymod{m_r} \cdots \gobymod{m_1} a_0)
=
(a^n_r \gobymod{m^n_r} \cdots \gobymod{m^n_1} a^n_0).
\]
Hence 
\[
t, t' \in \bigcap_{r \in \nat} \psi_{m_1} \cdots \psi_{m_r} (Ja_r),
\]
which by~\bref{item:pr-notmet} implies that $t = t'$, as required.

\paragraph*{}
The next part of the strategy is to show that $(J, \gamma)$ is the
universal solution in $\Set$.  Let $(X, \xi)$ be an $M$-coalgebra in
$\Set$.  Define for each $a \in \scat{A}$ and $x \in Xa$ a sequence
$(K_n(x))_{n\in\nat}$ of subsets of $Ja$ as follows:
\begin{itemize}
\item $K_0 (x) = Ja$
\item $K_{n+1} (x)
=
\bigcap 
\{ 
\psi_m K_n (y)
\such
b \gobymod{m} a, 
\ 
y \in Xb, 
\ 
\xi(x) = m \otimes y
\}$.
\end{itemize}
We eventually show that $\bigcap_{n\in\nat} K_n(x)$ is a one-element set
and that its element is the value at $x$ of the unique map $(X, \xi) \go
(J, \gamma)$.

\paragraph*{Claim}
$K_n (x) \supseteq K_{n+1}(x)$ for each $a \in \scat{A}$, $x \in Xa$, $n
\in \nat$.

\paragraph*{Proof}
For $n = 0$ this is trivial.  For $n\geq 0$, inductively,
\begin{eqnarray*}
K_{n+1}(x)	&
=	&
\bigcap 
\{
\psi_m K_n (y)
\such
\xi(x) = m \otimes y
\}	\\
	&\supseteq	&
\bigcap
\{
\psi_m K_{n+1}(y)
\such
\xi(x) = m \otimes y
\}	\\
	&=	&
K_{n+2}(x).
\end{eqnarray*}

\paragraph*{Claim}
$(Jf) K_n (x) \sub K_n (fx)$ for each $(a \goby{f} a')$ in $\scat{A}$, $x
\in Xa$, $n\in\nat$.

\paragraph*{Proof}
For $n = 0$ this is trivial.  Suppose inductively that the claim holds for
$n\in\nat$.  Suppose that $t \in K_{n+1}(x)$ and that $b' \gobymod{m'} a'$
and $y' \in Xb$ with $\xi (fx) = m' \otimes y'$; we have to show that $ft
\in \psi_{m'} K_n(y')$.

We may choose $b \gobymod{m} a$ and $y \in Xb$ such that $\xi(x) = m
\otimes y$.  Then $m' \otimes y' = \xi(fx) = fm \otimes y$, so there exist a
commutative square
\[
\begin{diagdiag}
	&		&c	&		&	\\
	&\ldTo<g	&	&\rdTo>{g'}	&	\\
b	&		&	&		&b'	\\
	&\rdMod<{fm}	&	&\ldMod>{m'}	&	\\
	&		&a'	&		&	\\
\end{diagdiag}
\]
and $z \in Xc$ such that $y = gz$ and $y' = g'z$.  Now $\xi(x) = m \otimes
y = m \otimes gz = mg \otimes z$, so $t \in \psi_{mg} K_n(z)$, so
\begin{eqnarray*}
ft	&
\in	&
\psi_{fmg} K_n(z)	
=
\psi_{m'g'} K_n(z)
=
\psi_{m'} (Jg') K_n(z)	\\
	&
\sub	&
\psi_{m'} K_n(g'z)
=	
\psi_{m'} K_n(y')
\end{eqnarray*}
(the penultimate step by inductive hypothesis), as required.

\paragraph*{Claim} 
$K_n(x)$ is closed in $Ja$ for each $a \in \scat{A}$, $x \in Xa$, $n \in
\nat$.

\paragraph*{Proof}
Lemma~\ref{lemma:fixedpointcomponents} and induction.

\paragraph*{Claim} 
$K_n(x)$ is nonempty for each $a \in \scat{A}$, $x \in Xa$, $n \in \nat$.

\paragraph*{Proof}
We use induction on $n\in\nat$ for all $a \in \scat{A}$ and $x \in Xa$
simultaneously.  

For $n = 0$, let $a \in \scat{A}$ and $x \in Xa$.  There exists a
resolution of $x$, and in particular an element of $\catI a$.  Since $J$ is
occupied, $\emptyset \neq Ja = K_0(x)$.

Now let $n\in\nat$, $a \in \scat{A}$, and $x \in Xa$; we have to prove that
$K_{n+1}(x)$ is nonempty.  Since $K_{n + 1}(x)$ is an intersection of a
family of closed subsets of a compact space, it suffices to show that the
intersection of any finite sub-family is nonempty.  So, suppose that $r \in
\nat$ and $\xi(x) = m_1 \otimes y_1 = \cdots = m_r \otimes y_r$ where $b_i
\gobymod{m_i} a$ and $y_i \in Xb_i$.  By
\cite[\lemmaequalityinMtensorX]{SS1} and an easy induction on $r$, there
exist $c \gobymod{p} a$ and $z \in Xc$ such that $\xi(x) = p \otimes z$,
and, for each $i \in \{ 1, \ldots, r \}$, a map $g_i: c \go b_i$ such that
$m_i g_i = p$ and $g_i z = y_i$.  But then for each $i$,
\begin{eqnarray*}
\psi_p K_n(z)	&
=	&
\psi_{m_i g_i} K_n(z) 
=
\psi_{m_i} (Jg_i) K_n(z)	\\
	&
\sub	&
\psi_{m_i} K_n (g_i z)
=
\psi_{m_i} K_n(y_i),
\end{eqnarray*}
so
\[
\psi_p K_n(z)
\sub
\psi_{m_1} K_n(y_1)
\cap \cdots \cap
\psi_{m_r} K_n(y_r),
\]
and $\psi_p K_n(z) \neq \emptyset$ by inductive hypothesis.

\paragraph*{Claim}
$\left| \bigcap_{n\in\nat} K_n(x) \right| = 1$
for each $a \in \scat{A}$, $x \in Xa$.

\paragraph*{Proof}
Let $a \in \scat{A}$ and $x \in Xa$.  The subsets $K_n(x)$ of the compact
space $Ja$ are nonempty, closed, and form a nested sequence, so
$\bigcap_{n\in\nat} K_n(x)$ has at least one element.  On the other hand,
we may choose a resolution
\[
( \cdots \gobymod{m_2} a_1 \gobymod{m_1} a_0, (x_n)_{n\in\nat})
\]
of $x$; then for each $n\in\nat$ we have $K_n(x) \sub \psi_{m_1} \cdots
\psi_{m_n} (Ja_n)$, so by~\bref{item:pr-notmet}, $\bigcap_{n\in\nat}
K_n(x)$ has at most one element.

\paragraph*{}%
By the last claim, there is for each $a \in \scat{A}$ a function
$\ovln{\xi}_a: Xa \go Ja$ defined by $\bigcap_{n\in\nat} K_n (x) = \{
\ovln{\xi}_a (x) \}$.

\paragraph*{Claim}
$(\ovln{\xi}_a)_{a \in \scat{A}}$ is a natural transformation $X \go J$ of
$\Set$-valued functors.

\paragraph*{Proof}
Let $f: a \go a'$ be a map in $\scat{A}$ and $x \in Xa$.  Then for all
$n\in\nat$,
\[
f \ovln{\xi}_a (x) 
\in
(Jf) K_n(x)
\sub
K_n (fx),
\]
so $f \ovln{\xi}_a (x) = \ovln{\xi}_{a'} (fx)$, as required.

\paragraph*{Claim}
$\ovln{\xi}: (X, \xi) \go (J, \gamma)$ is a map of coalgebras in $\Set$.

\paragraph*{Proof}
We have to show that for all $a \in \scat{A}$, the square
\[
\begin{diagram}
Xa			&\rTo^{\xi_a}		&(M \otimes X)a	\\
\dTo<{\ovln{\xi}_a}	&			&
\dTo>{(M \otimes \ovln{\xi})_a}					\\
Ja			&\rTo_{\gamma_a}	&(M \otimes J)a	\\
\end{diagram}
\]
commutes.  So let $x \in Xa$ and write $\xi_a(x) = (b \gobymod{m} a)
\otimes y$; we have to show that
$
\gamma_a \ovln{\xi}_a (x)
=
m \otimes \ovln{\xi}_b (y)
$,
or equivalently, 
$
\ovln{\xi}_a(x)
=
\psi_m \ovln{\xi}_b (y)
$,
or equivalently, 
$ 
\psi_m \ovln{\xi}_b (y)
\in 
K_n (x)
$
for all $n\in\nat$.
When $n = 0$ this is certainly true.  Now let $n\in\nat$; we have to show
that if $\xi_a(x) = (b' \gobymod{m'} a) \otimes y'$ then $\psi_m
\ovln{\xi}_b (y) \in \psi_{m'} K_n (y')$.  Since $m \otimes y = m' \otimes
y'$, there exist a commutative square
\[
\begin{diagdiag}
	&		&c	&		&	\\
	&\ldTo<g	&	&\rdTo>{g'}	&	\\
b	&		&	&		&b'	\\
	&\rdMod<m	&	&\ldMod>{m'}	&	\\
	&		&a	&		&	\\
\end{diagdiag}
\]
and $z \in Xc$ such that $gz = y$ and $g'z = y'$.  Hence
\begin{eqnarray*}
\psi_m \ovln{\xi}_b (y)	&
=	&
\psi_m \ovln{\xi}_b (gz)
=
\psi_m g \ovln{\xi}_c (z)	\\
	&
=	&
\psi_{mg} \ovln{\xi}_c (z)	
=	
\psi_{m'g'} \ovln{\xi}_c (z)	\\
	&
=	&
\psi_{m'} \ovln{\xi}_{b'} (y')
\in	
\psi_{m'} K_n(y')
\end{eqnarray*}
(the last equality by symmetry), as required.

\paragraph*{Claim} 
If $\twid{\xi}: (X, \xi) \go (J, \gamma)$ is a map of coalgebras in $\Set$
then $\twid{\xi} = \ovln{\xi}$.

\paragraph*{Proof}
We show by induction on $n\in\nat$ that $\twid{\xi}_a (x) \in K_n(x)$ for
all $a \in \scat{A}$ and $x \in Xa$; the result follows.  For $n=0$ this is
trivial.  Let $n \in \nat$, $a \in \scat{A}$, and $x \in Xa$.  If $\xi_a(x)
= (b \gobymod{m} a) \otimes y$ then, as in the proof of the previous claim,
$\twid{\xi}$ being a map of coalgebras implies that $\twid{\xi}_a (x) =
\psi_m \twid{\xi}_b (y)$, so by inductive hypothesis $\twid{\xi}_a (x) \in
\psi_m K_n(y)$.  Hence $\twid{\xi}_a (x) \in K_{n+1}(x)$, as required.

\paragraph*{}%
So $(J, \gamma)$ (or properly, $U_*(J, \gamma)$) is the terminal coalgebra
in $\Set$.  By Theorem~\thmexistenceofuniversalsolution\ of \cite{SS1},
condition \So\ holds, so $U_*(J, \gamma)$ is the universal solution $U_*(I,
\iota)$ constructed in~\S\secunivsoln\ of \cite{SS1}.  There is a
unique map $(J, \gamma) \go (I, \iota)$ of coalgebras in $\Top$, each
component $Ja \go Ia$ of which is a continuous bijection of compact
Hausdorff spaces, hence a homeomorphism.  So $(J, \gamma)$ is isomorphic to
$(I, \iota)$ as a coalgebra in $\Top$, and is therefore the universal
solution in $\Top$.  \done
\end{proof}

In many examples the universal solution is especially easy to recognize.
\begin{cor}[Crude Recognition Theorem]
\label{cor:cruderecognition}
Suppose that $\scat{A}$ is finite.  Let $(J, \gamma = \psi^{-1})$ be a
fixed point of $M$ in $\Top$ such that for each $a \in \scat{A}$, $Ja$ is
nonempty and compact, and suppose that the spaces $Ja$ can be metrized in
such a way that for each $b \gobymod{m} a$, the map $Jb \goby{\psi_m} Ja$
is a contraction.  Then $(J, \gamma)$ is the universal solution of
$(\scat{A}, M)$.
\end{cor}

\begin{proof}
Since $\scat{A}$ and $M$ are finite, there are only finitely many module
elements $m$, so we may choose $\lambda < 1$ such that each $\psi_m$ is a
contraction with constant $\lambda$.  Since $\scat{A}$ is finite and each
$Ja$ compact, we may also choose $D \geq 0$ such that $\diam(Ja) \leq D$
for all $a \in \scat{A}$.

We verify condition~\bref{item:pr-met} of the Precise Recognition Theorem.
Certainly $J$ is occupied.  For the main part of the condition, let $\epsln
> 0$; then we may choose $n \in \nat$ such that $\lambda^n D \leq \epsln$,
and the result follows.  
\done
\end{proof}

\section{Examples}
\label{sec:examples}

\begin{example*}{Interval}	\label{eg:rec-Freyd}
We finally prove the Topological Freyd
Theorem~\cite[\thmtopologicalFreyd]{SS1}.  So far we have verified that the
$(\scat{A}, M)$ concerned is a self-similarity system and exhibited an
$M$-coalgebra $(J, \gamma)$ (previously written $(I, \iota)$) with $J(1) =
[0, 1]$.  Evidently $\gamma$ is invertible, so we have a fixed point $(J,
\gamma = \psi^{-1})$.  The only non-constant maps $\psi_m$ are
\[
\begin{array}{ccc}
[0, 1]	&\go		&[0, 1]		\\
t	&\goesto	&t/2		\\
t	&\goesto	&(t + 1)/2,	
\end{array}
\]
both of which are contractions, so by the Crude Recognition Theorem, $(J,
\gamma)$ is the universal solution.

Freyd's Theorem expresses $[0, 1]$ as two copies of itself glued end to
end.  Two can be replaced by any larger number.  Thus, for each $k \geq 2$
there is a self-similarity structure $M^{(k)}$ on $\scat{A}$ (with, for
instance, $|M^{(k)} (1, 1)| = k$), the multiplication map $k \cdot \dashbk:
[0, 1] \go [0, k]$ puts an $M^{(k)}$-coalgebra structure $\gamma^{(k)}$ on
$J$, and the same argument shows that $(J, \gamma^{(k)})$ is the universal
solution. 
\end{example*}

\begin{example*}{Circle}
In most of our examples, every arrow in $\scat{A}$ is monic.  Here is an
exception.  Let $\scat{A}$ be the category generated by objects and arrows
\[
0 \parpair{\sigma}{\tau} 1 \goby{\rho} 2
\]
subject to $\rho \sigma = \rho \tau$, as in \cite[\egnondegencofork]{SS1}.
There is a self-similarity structure $M$ on $\scat{A}$ given informally by
the equations
\begin{eqnarray*}
X_0	&=	&X_0	\\
\begin{array}{c}
\setlength{\unitlength}{1mm}
\begin{picture}(30,14)(-15,-5)
\cell{0}{0}{c}{\includegraphics{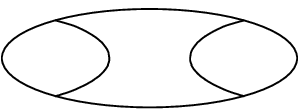}}
\cell{-10}{0}{c}{X_0}
\cell{10}{0}{c}{X_0}
\cell{0}{5.5}{b}{X_1}
\end{picture}
\end{array}
&
=	&
\begin{array}{c}
\setlength{\unitlength}{1mm}
\begin{picture}(49,14)(-24.5,-5)
\cell{0}{0}{c}{\includegraphics{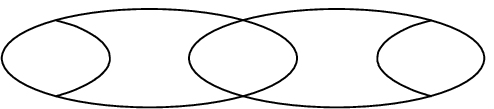}}
\cell{0}{0}{c}{X_0}
\cell{-20}{0}{c}{X_0}
\cell{20}{0}{c}{X_0}
\cell{-10}{5.5}{b}{X_1}
\cell{10}{5.5}{b}{X_1}
\end{picture}
\end{array}
\\
X_2	&=	&
\begin{array}{c}
\setlength{\unitlength}{1mm}
\begin{picture}(32,23)(-16,0)
\cell{0}{0}{b}{\includegraphics{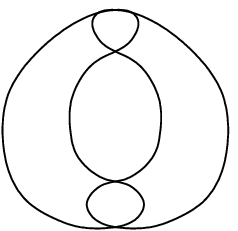}}
\cell{0}{2.5}{c}{X_0}
\cell{0}{20.5}{c}{X_0}
\cell{-11.5}{7}{r}{X_1}
\cell{11.5}{7}{l}{X_1}
\end{picture}
\end{array}
\end{eqnarray*}
and formally by
\[
\begin{diagram}
				&\ \ \				&
M(\dashbk, 0)			&\pile{\rTo^{\sigma\cdot\dashbk}\\ 
				\rTo_{\tau\cdot\dashbk}}	&
M(\dashbk, 1)			&\rTo^{\rho\cdot\dashbk}	&
M(\dashbk, 2)			\\
				&\				&
				&				&
				&				&
				\\
M(0, \dashbk)			&				&
\{ \id \}			&\pile{\rTo^0\\ \rTo_1}		&
\{ 0, \dhalf, 1 \}		&\rTo^{\mr{reduce\ } \mod 1}&
\{ 0, \dhalf \}			\\
\uTo<{\dashbk\cdot\sigma} 
\uTo>{\dashbk\cdot\tau}		&				&
\uTo \uTo			&				&
\uTo<\inf \uTo>\sup		&				&
\uTo<\inf \uTo>{\sup\ (\mod 1)}\\
M(1, \dashbk)			&				&
\emptyset			&\pile{\rTo\\ \rTo}		&
\{ [0, \dhalf], [\dhalf, 1] \}	&\rTo^{\id}			&
\{ [0, \dhalf], [\dhalf, 1] \}	\\
\uTo<{\dashbk\cdot\rho}		&				&
\uTo				&				&
\uTo				&				&
\uTo				\\
M(2, \dashbk)			&				&
\emptyset			&\pile{\rTo\\ \rTo}		&
\emptyset			&\rTo				&
\emptyset.			\\	
\end{diagram}
\]
(The top-left quarter is the Freyd module: eq.~(\eqFreydSSS)
of~\cite{SS1}.)  By \cite[\egnondegencofork]{SS1}, $M$ is nondegenerate,
and clearly $M$ is finite, so $(\scat{A}, M)$ is a self-similarity system.
We also have a nondegenerate functor $J: \scat{A} \go \Top$ given by
\[
\{\star\} \parpair{0}{1} [0, 1] \goby{q} S^1
\]
where $q$ identifies the endpoints.  Now $M \otimes J$ can naturally be
identified with
\[
\{\star\} \parpair{0}{2} [0, 2] \goby{\mr{quotient}} [0, 2] / (0 = 2),
\]
so there is an obvious isomorphism $\gamma: J \goiso M \otimes J$.
Moreover, each of the spaces $Ja$ is compact and nonempty.  The non-trivial
maps of the form $\psi_m$ (where $\psi = \gamma^{-1}$) are the two maps
\[
\psi_{[0, \half]},
\,
 \psi_{[\half, 1]}: 
[0, 1] 
\go
[0, 1]
\]
seen in Example~\ref{eg:rec-Freyd} and their respective composites with $q:
[0, 1] \go S^1$.  If we give $[0, 1]$ the standard metric and $S^1$ the
standard metric scaled down by a sufficiently large factor then all four
maps are contractions, so by the Crude Recognition Theorem, $(J, \gamma)$
is the universal solution.
\end{example*}

\begin{example*}{Cantor set}	\label{eg:rec-terminal}
Write $\One$ for the terminal category (one object and only the identity
arrow).  A self-similarity structure on $\One$ is a finite set $M$, and an
$M$-coalgebra is a space $X$ equipped with a map into the $M$-fold
coproduct $M \times X$.  The universal solution is the power $M^\nat$
(regarding $M$ as a discrete space) together with the shift isomorphism
$\gamma = \psi^{-1}: M^\nat \goiso M \times M^\nat$.  This can be seen
directly from the construction of the universal solution in
\cite[\S\secunivsoln]{SS1}, or from a Recognition Theorem as follows.
The space $M^\nat$ is nonempty if $M$ is, hence occupied.  It is compact.
For $m \in M$, the map $\psi_m: M^\nat \go M^\nat$ is given by
\[
\psi_m (m_0, m_1, \ldots)
=
(m, m_0, m_1, \ldots),
\]
so condition~\bref{item:pr-notmet} of the Precise Recognition Theorem
holds, as required.  When $M$ has cardinality $2$, the universal solution
is the standard Cantor set $\{0, 1\}^\nat$, often regarded as a subspace of
$[0, 1]$ via the embedding $(m_n)_{n\in\nat} \goesto \sum_{n\in\nat} 2 m_n
\cdot 3^{-(n+1)}$.

In fact, the homeomorphism type of $M^\nat$ is independent of $M$ for $|M|
\geq 2$.  This can be proved directly or as follows.  Write $k = \{ 0,
\ldots, k-1 \}$, write $\psi: 2 \times 2^\nat \goiso 2^\nat$ for the
(inverse) shift isomorphism, and define for each $k \geq 1$ an isomorphism
$\psi^{(k)}: k \times 2^\nat \goiso 2^\nat$ by
\begin{itemize}
\item $\psi^{(1)}$ is the canonical isomorphism
\item $\psi^{(k + 1)} 
=
\left(
(k + 1) \times 2^\nat
\rTo_\diso^{\mr{canonical}}
(k \times 2^\nat) + 2^\nat
\rTo_\diso^{\psi^{(k)} + 1}
2^\nat + 2^\nat
\rTo_\diso^{\psi}
2^\nat
\right). 
$
\end{itemize}
Then for each $k\geq 1$ and $m \in k$, the map $\psi^{(k)}_m: 2^\nat \go
2^\nat$ is of the form $\psi_{p_1} \cdots \psi_{p_r}$ with $r \geq 0$ and
$p_j \in 2$; if $k \geq 2$ then $r \geq 1$.  Using the metric on $2^\nat$
induced by its embedding into $[0, 1]$, $\psi_0$ and $\psi_1$ are
contractions with constant $1/3$, so if $k \geq 2$ then each $\psi^{(k)}_m$
is also a contraction with constant $1/3$.  So by the Crude Recognition
Theorem, $( 2^\nat, (\psi^{(k)})^{-1} )$ is the universal solution of
$(\One, k)$ whenever $k \geq 2$.  In particular, $2^\nat \iso k^\nat$.
\end{example*}

\begin{example*}{Discrete systems}	\label{eg:rec-discrete}
In~\S\ref{sec:thecantorset} we will look systematically at self-similarity
systems $(\scat{A}, M)$ in which $\scat{A}$ is a discrete category.
Proposition~\ref{propn:emptyorcantor} says that the universal solution
often consists simply of copies of the Cantor set and $\emptyset$.  But
it is not always so: for instance, the universal solution of the system
\begin{eqnarray*}
A	&=	&A	\\
B	&=	&A + B	
\end{eqnarray*}
(defined formally before eq.~(\eqinformalsss) of~\cite{SS1}) is $A = \{ 0
\}$ and $B = \nat \cup \{ \infty \}$, with $B$ topologized as the
Alexandroff compactification of the discrete space $\nat$.
\end{example*}

Finite products of self-similarity systems can be formed in the expected
way.
\begin{lemma}
\label{lemma:functorsonproducts}
Let $\left( \scat{B}_\omega \goby{Z_\omega} \Set \right)_{\omega \in
\Omega}$ be a family of functors on small categories $\scat{B}_\omega$,
indexed over an arbitrary set $\Omega$.  Then there is a functor
$\prod_{\omega \in \Omega} Z_\omega: \prod_{\omega \in \Omega}
\scat{B}_\omega \go \Set$, and
\begin{enumerate}
\item	\label{item:prod-elts}
$\elt{\prod Z_\omega} \iso \prod \elt{Z_\omega}$
\item \label{item:prod-fin} 
if $\Omega$ is finite and each $\elt{Z_\omega}$ is finite then so is
$\elt{\prod Z_\omega}$
\item	\label{item:prod-nondegen}
if each $Z_\omega$ is nondegenerate then so is $\prod Z_\omega$.
\end{enumerate}
\end{lemma}
\begin{proof}
\bref{item:prod-elts} is straightforward, and~\bref{item:prod-fin}
and~\bref{item:prod-nondegen} follow.
\done
\end{proof}

\begin{lemma}
\label{lemma:productandtensor}
Let $ \left( \scat{B}_\omega \goby{X_\omega} \Top \right)_{\omega \in
\Omega} $ and $ \left( \scat{B}_\omega^\op \goby{Y_\omega} \Set
\right)_{\omega \in \Omega} $ be families of functors on small categories
$\scat{B}_\omega$, where $\Omega$ is finite and each $X_\omega b$ is
compact Hausdorff.  Then
\[
\left(
\prod_{\omega \in \Omega} Y_\omega
\right)
\otimes
\left(
\prod_{\omega \in \Omega} X_\omega
\right)
\iso
\prod_{\omega \in \Omega}
(Y_\omega \otimes X_\omega).
\]
\end{lemma}

\begin{proof}
This is clear when $\Omega = \emptyset$, so by induction, it suffices to
consider $\Omega = \{ 1, 2 \}$.  In that case, 
\begin{eqnarray*}
(Y_1 \times Y_2) \otimes (X_1 \times X_2)	&
\iso	&
\int^{(b_1, b_2)} 
Y_1 b_1 \times Y_2 b_2 \times X_1 b_1 \times X_2 b_2	\\
	&\iso	&
\int^{b_1} 
Y_1 b_1 \times X_1 b_1 \times 
\int^{b_2} 
Y_2 b_2 \times X_2 b_2					\\
	&\iso	&
(Y_1 \otimes X_1) \times (Y_2 \otimes X_2),
\end{eqnarray*}
using in the second isomorphism the fact that if $K$ is compact Hausdorff
then $K \times \dashbk: \Top \go \Top$ preserves colimits.
\done
\end{proof}

\begin{propn}[Product self-similarity system]
\label{propn:productsss}
Let $((\scat{A}_\lambda, M_\lambda))_{\lambda \in \Lambda}$ be a finite
family of self-similarity systems.  Then there is a product self-similarity
system $(\prod_{\lambda \in \Lambda} \scat{A}_\lambda, \prod_{\lambda \in
\Lambda} M_\lambda)$, where
\[
\left(
\prod_{\lambda \in \Lambda} M_\lambda
\right)
( (b_\lambda)_{\lambda \in \Lambda}, 
(a_\lambda)_{\lambda \in \Lambda})
=
\prod_{\lambda \in \Lambda}
M_\lambda (b_\lambda, a_\lambda).
\]
If $(\scat{A}_\lambda, M_\lambda)$ has a universal solution $(I_\lambda,
\iota_\lambda)$ for each $\lambda \in \Lambda$ then $(\prod \scat{A}_\lambda,
\prod M_\lambda)$ has a universal solution $(\prod I_\lambda, \prod
\iota_\lambda)$.  
\end{propn}

\begin{proof}
Lemma~\ref{lemma:functorsonproducts} implies that $(\prod
\scat{A}_\lambda, \prod M_\lambda)$ is a self-similarity system.  Now
suppose that each $(\scat{A}_\lambda, M_\lambda)$ has a universal solution
$(I_\lambda, \iota_\lambda)$ in $\Top$.  The functor $\prod I_\lambda:
\prod \scat{A}_\lambda \go \Set$ is nondegenerate, since
\[
U \of (\prod I_\lambda)
\iso
\prod (U \of I_\lambda):
\prod \scat{A}_\lambda
\go 
\Set
\]
is nondegenerate by
Lemma~\ref{lemma:functorsonproducts}\bref{item:prod-nondegen} and each
space $(\prod_{\lambda \in \Lambda} I_\lambda) (a_\lambda)_{\lambda \in
\Lambda}$ is compact Hausdorff.  By Lemma~\ref{lemma:productandtensor}, we
have a coalgebra structure
\[
\prod \iota_\lambda:
\prod I_\lambda
\go
\prod (M_\lambda \otimes I_\lambda)
\iso
(\prod M_\lambda) \otimes (\prod I_\lambda)
\]
on $\prod I_\lambda$, and $\prod \iota_\lambda$ is an isomorphism.  Also,
$\prod I_\lambda$ is occupied since each $I_\lambda$ is.  To finish the
proof it remains only to verify that $(\prod I_\lambda, \prod
\iota_\lambda)$ satisfies the main condition in~\bref{item:pr-notmet} of
the Precise Recognition Theorem, and this follows from the fact that it is
satisfied by each $(I_\lambda, \iota_\lambda)$.  \done
\end{proof}


\begin{example*}{Cubes}	\label{eg:rec-cubes}
Let $(\scat{A}, M)$ be the Freyd self-similarity system.  Then by
Proposition~\ref{propn:productsss}, the universal solution of $(\scat{A}^2,
M^2)$ has a universal solution $(I, \iota)$ satisfying $I(1, 1) = [0,
1]^2$.  Informally, the self-similarity equations are
\[
\begin{array}{rclcrcl}
\zmark		&=	&\zmark	&
\diagspace	&
\begin{array}{c}
\setlength{\unitlength}{1em}
\begin{picture}(3,0.5)(-1.5,-0.25)
\put(-1.3,0){\line(1,0){2.6}}
\cell{-1.3}{0}{c}{\zmark}
\cell{1.3}{0}{c}{\zmark}
\end{picture}
\end{array}
&
=	&
\begin{array}{c}
\setlength{\unitlength}{1em}
\begin{picture}(5.6,0.5)(-2.8,-0.25)
\put(0,0){\line(1,0){2.6}}
\put(0,0){\line(-1,0){2.6}}
\cell{-2.6}{0}{c}{\zmark}
\cell{0}{0}{c}{\zmark}
\cell{2.6}{0}{c}{\zmark}
\end{picture}
\end{array}
\\
&&&&&&\\
\begin{array}{c}
\setlength{\unitlength}{1em}
\begin{picture}(0.5,3)(-0.25,-1.5)
\put(0,-1.3){\line(0,1){2.6}}
\cell{0}{-1.3}{c}{\zmark}
\cell{0}{1.3}{c}{\zmark}
\end{picture}
\end{array}
&
=	&
\begin{array}{c}
\setlength{\unitlength}{1em}
\begin{picture}(0.5,5.6)(-0.25,-2.8)
\put(0,0){\line(0,1){2.6}}
\put(0,0){\line(0,-1){2.6}}
\cell{0}{-2.6}{c}{\zmark}
\cell{0}{0}{c}{\zmark}
\cell{0}{2.6}{c}{\zmark}
\end{picture}
\end{array}
&
\diagspace	&
\begin{array}{c}
\setlength{\unitlength}{1em}
\begin{picture}(3,3)(-1.5,-1.5)
\put(-1.3,-1.3){\textcolor{grey}{\rule{2.6\unitlength}{2.6\unitlength}}}
\put(-1.3,-1.3){\line(1,0){2.6}}
\put(-1.3,1.3){\line(1,0){2.6}}
\put(-1.3,-1.3){\line(0,1){2.6}}
\put(1.3,-1.3){\line(0,1){2.6}}
\cell{-1.3}{-1.3}{c}{\zmark}
\cell{1.3}{-1.3}{c}{\zmark}
\cell{-1.3}{1.3}{c}{\zmark}
\cell{1.3}{1.3}{c}{\zmark}
\end{picture}
\end{array}
&
=	&
\begin{array}{c}
\setlength{\unitlength}{1em}
\begin{picture}(5.6,5.6)(-2.8,-2.8)
\put(-2.6,-2.6){\textcolor{grey}{\rule{5.2\unitlength}{5.2\unitlength}}}
\put(-2.6,-2.6){\line(0,1){5.2}}
\put(0,-2.6){\line(0,1){5.2}}
\put(2.6,-2.6){\line(0,1){5.2}}
\put(-2.6,-2.6){\line(1,0){5.2}}
\put(-2.6,0){\line(1,0){5.2}}
\put(-2.6,2.6){\line(1,0){5.2}}
\cell{-2.6}{-2.6}{c}{\zmark}
\cell{0}{-2.6}{c}{\zmark}
\cell{2.6}{-2.6}{c}{\zmark}
\cell{-2.6}{0}{c}{\zmark}
\cell{0}{0}{c}{\zmark}
\cell{2.6}{0}{c}{\zmark}
\cell{-2.6}{2.6}{c}{\zmark}
\cell{0}{2.6}{c}{\zmark}
\cell{2.6}{2.6}{c}{\zmark}
\end{picture}
\end{array}
.
\end{array}
\]
For general $n\in\nat$, $(\scat{A}^n, M^n)$ has universal solution given by
$[0, 1]^n$.
\end{example*}

\begin{example}
Similarly, if $C$ denotes the Cantor set then the space
\[
C \times [0, 1]	
=	
\begin{array}{c}
\setlength{\unitlength}{1em}
\begin{picture}(2.7,2.9)(0,-0.1)
\put(0.1,0){\line(0,1){2.7}}
\put(0.2,0){\line(0,1){2.7}}
\put(0.7,0){\line(0,1){2.7}}
\put(0.8,0){\line(0,1){2.7}}
\put(1.9,0){\line(0,1){2.7}}
\put(2.0,0){\line(0,1){2.7}}
\put(2.5,0){\line(0,1){2.7}}
\put(2.6,0){\line(0,1){2.7}}
\put(0.3,0){\line(0,1){2.7}}
\put(0.6,0){\line(0,1){2.7}}
\put(2.1,0){\line(0,1){2.7}}
\put(2.4,0){\line(0,1){2.7}}
\put(0.9,0){\line(0,1){2.7}}
\put(1.8,0){\line(0,1){2.7}}
\put(0,0){\line(0,1){2.7}}
\put(2.7,0){\line(0,1){2.7}}
\end{picture}
\end{array}
\]
arises from the self-similarity system
\[
\begin{array}{rclcrcl}
\cdot\	&
=	&
\ \cdot	&
\diagspace	&
\begin{array}{c}	
\setlength{\unitlength}{1em}
\begin{picture}(2.7,0.2)(0,-0.1)
\cell{0.1}{0}{c}{\cdot}
\cell{0.2}{0}{c}{\cdot}
\cell{0.7}{0}{c}{\cdot}
\cell{0.8}{0}{c}{\cdot}
\cell{1.9}{0}{c}{\cdot}
\cell{2.0}{0}{c}{\cdot}
\cell{2.5}{0}{c}{\cdot}
\cell{2.6}{0}{c}{\cdot}
\cell{0.3}{0}{c}{\cdot}
\cell{0.6}{0}{c}{\cdot}
\cell{2.1}{0}{c}{\cdot}
\cell{2.4}{0}{c}{\cdot}
\cell{0.9}{0}{c}{\cdot}
\cell{1.8}{0}{c}{\cdot}
\cell{0}{0}{c}{\cdot}
\cell{2.7}{0}{c}{\cdot}
\end{picture}
\end{array}
&
=	&
\begin{array}{c}
\setlength{\unitlength}{1em}
\begin{picture}(8.1,0.2)(0,-0.1)
\cell{0.1}{0}{c}{\cdot}
\cell{0.2}{0}{c}{\cdot}
\cell{0.7}{0}{c}{\cdot}
\cell{0.8}{0}{c}{\cdot}
\cell{1.9}{0}{c}{\cdot}
\cell{2.0}{0}{c}{\cdot}
\cell{2.5}{0}{c}{\cdot}
\cell{2.6}{0}{c}{\cdot}
\cell{0.3}{0}{c}{\cdot}
\cell{0.6}{0}{c}{\cdot}
\cell{2.1}{0}{c}{\cdot}
\cell{2.4}{0}{c}{\cdot}
\cell{0.9}{0}{c}{\cdot}
\cell{1.8}{0}{c}{\cdot}
\cell{0}{0}{c}{\cdot}
\cell{2.7}{0}{c}{\cdot}
\cell{5.5}{0}{c}{\cdot}
\cell{5.6}{0}{c}{\cdot}
\cell{6.1}{0}{c}{\cdot}
\cell{6.2}{0}{c}{\cdot}
\cell{7.3}{0}{c}{\cdot}
\cell{7.4}{0}{c}{\cdot}
\cell{7.9}{0}{c}{\cdot}
\cell{8.0}{0}{c}{\cdot}
\cell{5.7}{0}{c}{\cdot}
\cell{6.0}{0}{c}{\cdot}
\cell{7.5}{0}{c}{\cdot}
\cell{7.8}{0}{c}{\cdot}
\cell{6.3}{0}{c}{\cdot}
\cell{7.2}{0}{c}{\cdot}
\cell{5.4}{0}{c}{\cdot}
\cell{8.1}{0}{c}{\cdot}
\end{picture}
\end{array}
\\
&&&&&&\\
\begin{array}{c}
\setlength{\unitlength}{1em}
\begin{picture}(0,2.7)(0,0)
\put(0,0){\line(0,1){2.7}}
\cell{0}{0}{c}{\cdot}
\cell{0}{2.7}{c}{\cdot}
\end{picture}
\end{array}
&
=	&
\begin{array}{c}
\setlength{\unitlength}{1em}
\begin{picture}(0,5.4)(0,0)
\put(0,0){\line(0,1){5.4}}
\cell{0}{0}{c}{\cdot}
\cell{0}{2.7}{c}{\cdot}
\cell{0}{5.4}{c}{\cdot}
\end{picture}
\end{array}
&
\diagspace	&
\begin{array}{c}	
\setlength{\unitlength}{1em}
\begin{picture}(2.7,2.9)(0,-0.1)
\put(0.1,0){\line(0,1){2.7}}
\put(0.2,0){\line(0,1){2.7}}
\put(0.7,0){\line(0,1){2.7}}
\put(0.8,0){\line(0,1){2.7}}
\put(1.9,0){\line(0,1){2.7}}
\put(2.0,0){\line(0,1){2.7}}
\put(2.5,0){\line(0,1){2.7}}
\put(2.6,0){\line(0,1){2.7}}
\put(0.3,0){\line(0,1){2.7}}
\put(0.6,0){\line(0,1){2.7}}
\put(2.1,0){\line(0,1){2.7}}
\put(2.4,0){\line(0,1){2.7}}
\put(0.9,0){\line(0,1){2.7}}
\put(1.8,0){\line(0,1){2.7}}
\put(0,0){\line(0,1){2.7}}
\put(2.7,0){\line(0,1){2.7}}
\cell{0.1}{0}{c}{\cdot}
\cell{0.2}{0}{c}{\cdot}
\cell{0.7}{0}{c}{\cdot}
\cell{0.8}{0}{c}{\cdot}
\cell{1.9}{0}{c}{\cdot}
\cell{2.0}{0}{c}{\cdot}
\cell{2.5}{0}{c}{\cdot}
\cell{2.6}{0}{c}{\cdot}
\cell{0.3}{0}{c}{\cdot}
\cell{0.6}{0}{c}{\cdot}
\cell{2.1}{0}{c}{\cdot}
\cell{2.4}{0}{c}{\cdot}
\cell{0.9}{0}{c}{\cdot}
\cell{1.8}{0}{c}{\cdot}
\cell{0}{0}{c}{\cdot}
\cell{2.7}{0}{c}{\cdot}
\cell{0.1}{2.7}{c}{\cdot}
\cell{0.2}{2.7}{c}{\cdot}
\cell{0.7}{2.7}{c}{\cdot}
\cell{0.8}{2.7}{c}{\cdot}
\cell{1.9}{2.7}{c}{\cdot}
\cell{2.0}{2.7}{c}{\cdot}
\cell{2.5}{2.7}{c}{\cdot}
\cell{2.6}{2.7}{c}{\cdot}
\cell{0.3}{2.7}{c}{\cdot}
\cell{0.6}{2.7}{c}{\cdot}
\cell{2.1}{2.7}{c}{\cdot}
\cell{2.4}{2.7}{c}{\cdot}
\cell{0.9}{2.7}{c}{\cdot}
\cell{1.8}{2.7}{c}{\cdot}
\cell{0}{2.7}{c}{\cdot}
\cell{2.7}{2.7}{c}{\cdot}
\end{picture}
\end{array}
&
=	&
\begin{array}{c}
\setlength{\unitlength}{1em}
\begin{picture}(8.1,5.6)(0,-0.1)
\put(0.1,0){\line(0,1){5.4}}
\put(0.2,0){\line(0,1){5.4}}
\put(0.7,0){\line(0,1){5.4}}
\put(0.8,0){\line(0,1){5.4}}
\put(1.9,0){\line(0,1){5.4}}
\put(2.0,0){\line(0,1){5.4}}
\put(2.5,0){\line(0,1){5.4}}
\put(2.6,0){\line(0,1){5.4}}
\put(0.3,0){\line(0,1){5.4}}
\put(0.6,0){\line(0,1){5.4}}
\put(2.1,0){\line(0,1){5.4}}
\put(2.4,0){\line(0,1){5.4}}
\put(0.9,0){\line(0,1){5.4}}
\put(1.8,0){\line(0,1){5.4}}
\put(0,0){\line(0,1){5.4}}
\put(2.7,0){\line(0,1){5.4}}
\put(5.5,0){\line(0,1){5.4}}
\put(5.6,0){\line(0,1){5.4}}
\put(6.1,0){\line(0,1){5.4}}
\put(6.2,0){\line(0,1){5.4}}
\put(7.3,0){\line(0,1){5.4}}
\put(7.4,0){\line(0,1){5.4}}
\put(7.9,0){\line(0,1){5.4}}
\put(8.0,0){\line(0,1){5.4}}
\put(5.7,0){\line(0,1){5.4}}
\put(6.0,0){\line(0,1){5.4}}
\put(7.5,0){\line(0,1){5.4}}
\put(7.8,0){\line(0,1){5.4}}
\put(6.3,0){\line(0,1){5.4}}
\put(7.2,0){\line(0,1){5.4}}
\put(5.4,0){\line(0,1){5.4}}
\put(8.1,0){\line(0,1){5.4}}
\cell{0.1}{0}{c}{\cdot}
\cell{0.2}{0}{c}{\cdot}
\cell{0.7}{0}{c}{\cdot}
\cell{0.8}{0}{c}{\cdot}
\cell{1.9}{0}{c}{\cdot}
\cell{2.0}{0}{c}{\cdot}
\cell{2.5}{0}{c}{\cdot}
\cell{2.6}{0}{c}{\cdot}
\cell{0.3}{0}{c}{\cdot}
\cell{0.6}{0}{c}{\cdot}
\cell{2.1}{0}{c}{\cdot}
\cell{2.4}{0}{c}{\cdot}
\cell{0.9}{0}{c}{\cdot}
\cell{1.8}{0}{c}{\cdot}
\cell{0}{0}{c}{\cdot}
\cell{2.7}{0}{c}{\cdot}
\cell{5.5}{0}{c}{\cdot}
\cell{5.6}{0}{c}{\cdot}
\cell{6.1}{0}{c}{\cdot}
\cell{6.2}{0}{c}{\cdot}
\cell{7.3}{0}{c}{\cdot}
\cell{7.4}{0}{c}{\cdot}
\cell{7.9}{0}{c}{\cdot}
\cell{8.0}{0}{c}{\cdot}
\cell{5.7}{0}{c}{\cdot}
\cell{6.0}{0}{c}{\cdot}
\cell{7.5}{0}{c}{\cdot}
\cell{7.8}{0}{c}{\cdot}
\cell{6.3}{0}{c}{\cdot}
\cell{7.2}{0}{c}{\cdot}
\cell{5.4}{0}{c}{\cdot}
\cell{8.1}{0}{c}{\cdot}
\cell{0.1}{2.7}{c}{\cdot}
\cell{0.2}{2.7}{c}{\cdot}
\cell{0.7}{2.7}{c}{\cdot}
\cell{0.8}{2.7}{c}{\cdot}
\cell{1.9}{2.7}{c}{\cdot}
\cell{2.0}{2.7}{c}{\cdot}
\cell{2.5}{2.7}{c}{\cdot}
\cell{2.6}{2.7}{c}{\cdot}
\cell{0.3}{2.7}{c}{\cdot}
\cell{0.6}{2.7}{c}{\cdot}
\cell{2.1}{2.7}{c}{\cdot}
\cell{2.4}{2.7}{c}{\cdot}
\cell{0.9}{2.7}{c}{\cdot}
\cell{1.8}{2.7}{c}{\cdot}
\cell{0}{2.7}{c}{\cdot}
\cell{2.7}{2.7}{c}{\cdot}
\cell{5.5}{2.7}{c}{\cdot}
\cell{5.6}{2.7}{c}{\cdot}
\cell{6.1}{2.7}{c}{\cdot}
\cell{6.2}{2.7}{c}{\cdot}
\cell{7.3}{2.7}{c}{\cdot}
\cell{7.4}{2.7}{c}{\cdot}
\cell{7.9}{2.7}{c}{\cdot}
\cell{8.0}{2.7}{c}{\cdot}
\cell{5.7}{2.7}{c}{\cdot}
\cell{6.0}{2.7}{c}{\cdot}
\cell{7.5}{2.7}{c}{\cdot}
\cell{7.8}{2.7}{c}{\cdot}
\cell{6.3}{2.7}{c}{\cdot}
\cell{7.2}{2.7}{c}{\cdot}
\cell{5.4}{2.7}{c}{\cdot}
\cell{8.1}{2.7}{c}{\cdot}
\cell{0.1}{5.4}{c}{\cdot}
\cell{0.2}{5.4}{c}{\cdot}
\cell{0.7}{5.4}{c}{\cdot}
\cell{0.8}{5.4}{c}{\cdot}
\cell{1.9}{5.4}{c}{\cdot}
\cell{2.0}{5.4}{c}{\cdot}
\cell{2.5}{5.4}{c}{\cdot}
\cell{2.6}{5.4}{c}{\cdot}
\cell{0.3}{5.4}{c}{\cdot}
\cell{0.6}{5.4}{c}{\cdot}
\cell{2.1}{5.4}{c}{\cdot}
\cell{2.4}{5.4}{c}{\cdot}
\cell{0.9}{5.4}{c}{\cdot}
\cell{1.8}{5.4}{c}{\cdot}
\cell{0}{5.4}{c}{\cdot}
\cell{2.7}{5.4}{c}{\cdot}
\cell{5.5}{5.4}{c}{\cdot}
\cell{5.6}{5.4}{c}{\cdot}
\cell{6.1}{5.4}{c}{\cdot}
\cell{6.2}{5.4}{c}{\cdot}
\cell{7.3}{5.4}{c}{\cdot}
\cell{7.4}{5.4}{c}{\cdot}
\cell{7.9}{5.4}{c}{\cdot}
\cell{8.0}{5.4}{c}{\cdot}
\cell{5.7}{5.4}{c}{\cdot}
\cell{6.0}{5.4}{c}{\cdot}
\cell{7.5}{5.4}{c}{\cdot}
\cell{7.8}{5.4}{c}{\cdot}
\cell{6.3}{5.4}{c}{\cdot}
\cell{7.2}{5.4}{c}{\cdot}
\cell{5.4}{5.4}{c}{\cdot}
\cell{8.1}{5.4}{c}{\cdot}
\end{picture}
\end{array}.
\end{array}
\]
\end{example}

Various familiar self-similar spaces arise as follows.  An \demph{iterated
function system} on $\reals^d$ is a family $\psi_0, \ldots, \psi_n$ ($n\geq
0$) of contractions $\reals^d \go \reals^d$.  A theorem of
Hutchinson~\cite{Hut} implies that there is a unique nonempty compact
subset $S$ of $\reals^d$ satisfying $S = \bigcup_{i = 0}^n \psi_i S$.

\begin{example*}{Sierpi\'nski simplices}	\label{eg:rec-Sierpinski}
Let $n\in\nat$ and let $s_0, \ldots, s_n$ be affinely independent points of
$\reals^n$.  For each $i \in \{ 0, \ldots, n \}$, write $\psi_i: \reals^n
\go \reals^n$ for the scaling with fixed point $s_i$ and scale factor
$\half$.  The \demph{Sierpi\'nski simplex} with vertices $s_0, \ldots, s_n$
is the unique nonempty compact subset $S$ of $\reals^n$ satisfying $S =
\bigcup_{i = 0}^n \psi_i S$.  When $n = 1$, this is the closed interval with
endpoints $s_0$ and $s_1$.  When $n = 2$, this is the usual Sierpi\'nski
triangle or gasket, which satisfies the equation
\[
\begin{array}{c}
\setlength{\unitlength}{1mm}
\begin{picture}(12,10)(-6,0)
\cell{0}{0}{b}{\includegraphics{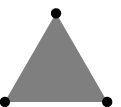}}
\cell{0}{4}{c}{S}
\end{picture}
\end{array}
=
\begin{array}{c}
\setlength{\unitlength}{1mm}
\begin{picture}(24,20)(-12,0)
\cell{0}{0}{b}{\includegraphics{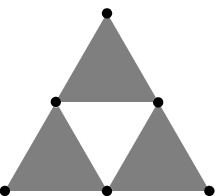}}
\cell{-5.5}{4}{c}{S}
\cell{5.5}{4}{c}{S}
\cell{0}{13}{c}{S}
\end{picture}
\end{array}.
\]

In arbitrary dimension, $S$ arises from a self-similarity system as
follows.  Let $\scat{A}$ be the category with objects $0$ and $1$ and
non-identity arrows $\sigma_0, \ldots, \sigma_n: 0 \go 1$.  Define an
equivalence relation $\sim$ on $\{ 0, \ldots, n \}$ by $(i, j) \sim (i',
j')$ if and only if $\{ i, j \} = \{ i', j' \}$, write $[i, j]$ for the
equivalence class of $(i, j)$, and define $M: \scat{A} \gomod \scat{A}$ by
\[
\begin{diagram}[height=3em,width=3em]
						&\ \ \		&
M(\dashbk, 0)							&
\pile{\rTo^{\sigma_0 \cdot \dashbk}\\ 
	\avdots\\ 
	\rTo_{\sigma_n \cdot \dashbk}}				&
M(\dashbk, 1)							\\
\						&		&
								&
								&
								\\
M(0, \dashbk)					&		&
\{ \id \}							&
\pile{\rTo^{[0, 0]}\\ \avdots\\ \rTo_{[n, n]}}			&
\{ 0, \ldots, n \}^2 / \sim					\\
\uTo<{\dashbk \cdot \sigma_0} 
\cdots 
\uTo>{\dashbk \cdot \sigma_n}			&		&
\uTo \cdots \uTo						&
								&
\uTo<{[\dashbk, 0]} \cdots \uTo>{[\dashbk, n]}			\\
M(1, \dashbk)					&		&
\emptyset							&
\pile{\rTo\\ \avdots\\ \rTo}					&
\{ 0, \ldots, n \}.						\\
\end{diagram}
\]
Then $(\scat{A}, M)$ is a self-similarity system.  Any space $X$ equipped
with distinct basepoints $x_0, \ldots, x_n$ defines a nondegenerate functor
$X: \scat{A} \go \Top$, and then $M \otimes X$ is the quotient space
\[
\frac{ 
\{ 0, \ldots, n \} \times X
}{
(i, x_j) = (j, x_i) \textrm{ for all } i, j
}
\]
with basepoints $(0, x_0), \ldots, (n, x_n)$.  Choose affinely independent
points $s_0, \ldots, s_n \in \reals^n$ and take $S$ as above.  Then $(S,
s_0, \ldots, s_n)$ determines a nondegenerate functor $J: \scat{A} \go
\Top$, and $\psi_0, \ldots, \psi_n$ determine a basepoint-preserving map
\begin{equation}	\label{eq:Sierpinski-map}
\begin{array}{ccc}
\{ 0, \ldots, n \} \times S	&\go		&S		\\
(i, s)				&\goesto	&\psi_i(s).
\end{array}
\end{equation}
If $i \neq j$ then the only point in $\psi_i S \cap \psi_j S$ is the
midpoint $\psi_i s_j = \psi_j s_i$ of $s_i$ and $s_j$.  Also, each $\psi_i$
is injective and, by definition of $S$, the map~\bref{eq:Sierpinski-map} is
surjective, so~\bref{eq:Sierpinski-map} induces a basepoint-preserving
homeomorphism $(M \otimes J) 1 \goiso J1$.  In other words, we have an
isomorphism $\psi: M \otimes J \goiso J$.  The spaces $J0 = \{\star\}$ and
$J1 = S$ are nonempty and compact and the structure maps $\psi_i: J1 \go
J1$ are contractions, so by the Crude Recognition Theorem, $(J, \psi^{-1})$
is the universal solution of $(\scat{A}, M)$.
\end{example*}

\begin{example*}{IFSs}	\label{eg:rec-ifs}
More generally, let $(\psi_0, \ldots, \psi_n)$ be an iterated function
system on $\reals^d$ (some $d\in\nat$); write $S$ for the unique nonempty
compact subset of $\reals^d$ satisfying $S = \bigcup_{i=0}^n \psi_i S$, and
$s_i$ for the fixed point of $\psi_i$.  Suppose that $\psi_0, \ldots,
\psi_n$ are injective, that $s_0, \ldots, s_n$ are distinct, and that if
$\psi_i s = \psi_j t$ with $i \neq j$ and $s, t \in S$ then $s, t \in \{
s_0, \ldots, s_n \}$.  Then the space $S$ arises from a self-similarity
system, as follows.

Define an equivalence relation $\sim$ on $\{ 0, \ldots, n \}^2$ by $(i, j)
\sim (i', j') \iff \psi_i s_j = \psi_{i'} s_{j'}$, and write $[i, j]$ for
the equivalence class of $(i, j)$.  Exactly as in the previous example,
this equivalence relation gives rise to a self-similarity system
$(\scat{A}, M)$, the space $S$ with basepoints $s_0, \ldots, s_n$
determines a nondegenerate functor $J: \scat{A} \go \Top$, the maps
$\psi_i$ determine an isomorphism $M \otimes J \goiso J$, and by the Crude
Recognition Theorem, this is the universal solution of $(\scat{A}, M)$.

Solutions to iterated function systems in which the small copies $\psi_i S$
of $S$ overlap by more than just a point can sometimes by described as
universal solutions too: $[0, 1]^n$ in Example~\ref{eg:rec-cubes}, for
instance.
\end{example*}

\begin{example*}{Barycentric subdivision}	\label{eg:rec-bary}
Barycentric subdivision expresses the $n$-simplex $\Delta^n$ as $(n + 1)!$
smaller copies of itself glued together along simplices of lower
dimension.  This self-similarity can be formalized as follows.

Let $\Dface$ be the category whose objects are the nonempty finite totally
ordered sets $\upr{n} = \{ 0, \ldots, n \}$ ($n \in \nat$) and whose maps are
the order-preserving injections.  For each $n, m \in \nat$, put
\[
M( \upr{n}, \upr{m} )
=
\{
\textrm{chains }
\emptyset \propersub Q(0) \propersub \cdots \propersub Q(n) \sub \upr{m}
\}
\]
where $\propersub$ means proper subset.  The idea can be seen in
Figure~\ref{fig:bary}:
\begin{figure}
\centering
\setlength{\unitlength}{1mm}
\begin{picture}(32,29)(-16,0)
\cell{0}{2.5}{b}{\includegraphics{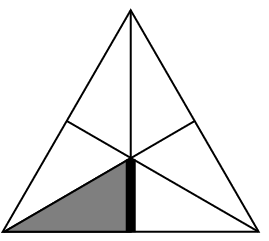}}
\cell{-16}{0}{bl}{0}
\cell{0}{29}{t}{1}
\cell{16}{0}{br}{2}
\cell{-8}{15}{br}{01}
\cell{8}{15}{bl}{12}
\cell{0}{0}{b}{02}
\cell{0}{12}{b}{012}
\end{picture}
\caption{Barycentric subdivision of $\Delta^2$}
\label{fig:bary}
\end{figure}
the $1$-simplex in bold and the shaded 2-simplex correspond respectively to
\begin{eqnarray*}
\left(
\emptyset \propersub \{ 0, 2 \} \propersub \{ 0, 1, 2 \}
\right)	&
\in	&
M(\upr{1}, \upr{2}),	\\
\left(
\emptyset \propersub \{ 0 \} \propersub \{ 0, 2 \} \propersub \{ 0, 1, 2 \}
\right)	&
\in	&
M(\upr{2}, \upr{2}).
\end{eqnarray*}
An element of $M(\upr{n}, \upr{m})$ can be regarded as an order-preserving
injection $\upr{n} \go \nepower \upr{m}$, where $\nepower$ denotes the set
of nonempty subsets ordered by inclusion.  By using direct images,
$\nepower \upr{m}$ is functorial in $\upr{m}$, so $M$ defines a module
$\Dface \gomod \Dface$.  It is straightforward to show that $\Dface$ has
all pullbacks and equalizers and that $\nepower: \Dface \go \fcat{Poset}$
preserves them, from which it follows that $M$ is nondegenerate
\cite[\corcategorywithpullbacksandequalizers]{SS1}.  And clearly $M$ is
finite, so $(\Dface, M)$ is a self-similarity system.

We are going to show that the universal solution is given by the simplex
functor $\Delta^\blob: \Dface \go \Top$.  To define $\Delta^\blob$, choose
for each $n\in\nat$ an affinely independent sequence $\vtr{e}^n_0, \ldots,
\vtr{e}^n_n$ of elements of $\reals^n$ and let $\Delta^n$ be their convex
hull; then for each map $f: \upr{n} \go \upr{m}$ in $\Dface$ there is a
unique affine map $\reals^n \go \reals^m$ sending $\vtr{e}^n_j$ to
$\vtr{e}^m_{f(j)}$ for each $j$, which restricts to a map $\Delta f = f_*:
\Delta^n \go \Delta^m$.  Again it is straightforward to check that $U \of
\Delta^\blob: \Dface \go \Set$ preserves pullbacks and equalizers, and each
space $\Delta^n$ is compact Hausdorff, so $\Delta^\blob$ is nondegenerate.

Given $Q: \upr{n} \gomod \upr{m}$, there is a unique affine map $\reals^n
\go \reals^m$ such that
\[
\vtr{e}^n_j
\goesto
\frac{1}{|Q(j)|}
\sum_{i \in Q(j)} \vtr{e}_i^m
\]
for all $j \in \upr{n}$, and this restricts to a map $\psi_Q: \Delta^n \go
\Delta^m$.  It can be checked that $\psi_{fQg} = f_* \of \psi_Q \of g_*$
for all $f$, $Q$, and $g$, so by the remarks at the beginning
of~\S\ref{sec:recognitiontheorems}, $\psi$ defines an $M$-algebra structure
on $\Delta^\blob$.

A standard calculation~\cite[2.21]{Hat} shows that for any
$Q: \upr{n} \gomod \upr{m}$,
\[
\diam (\psi_Q \Delta^n)
\leq
\frac{m}{m+1} \diam (\Delta^m)
\]
in the Euclidean metric.  More generally, if
\[
\upr{n_r} \gobymod{Q_r} \cdots \gobymod{Q_1} \upr{n_0}
\]
then the same method shows that
\begin{eqnarray*}
\diam (\psi_{Q_1} \cdots \psi_{Q_r} \Delta^{n_r})	&
\leq	&
\left( \frac{n_{r-1}}{n_{r-1} + 1} \right)
\cdots
\left( \frac{n_0}{n_0 + 1} \right)
\diam (\Delta^{n_0})	\\
	&
\leq	&
\left( \frac{n_0}{n_0 + 1} \right)^r
\diam (\Delta^{n_0}).
\end{eqnarray*}
So as long as $\psi$ is an isomorphism, condition~\bref{item:pr-met} of the
Precise Recognition Theorem holds and $(\Delta^\blob, \psi^{-1})$ is the
universal solution.

It remains to show that $\psi$ is an isomorphism.  Since the spaces
involved are compact Hausdorff, this says that for each $m \in \nat$,
\[
\psi: 
\int^{\upr{n}} M(\upr{n}, \upr{m}) \times \Delta^n
\go
\Delta^m
\]
is a bijection.  The proof is somewhat technical.

For surjectivity, take $\vtr{s} \in \Delta^m$; then $\vtr{s} = \sum_{i =
0}^m \sigma_i \vtr{e}^m_i$ with $\sigma_i \geq 0$ and $\sum \sigma_i = 1$.
There are unique $n \in \nat$ and $\kappa_0 > \cdots > \kappa_n >
\kappa_{n+1} = 0$ such that
\[
\{ \kappa_0, \ldots, \kappa_n, \kappa_{n+1} \}
=
\{ \sigma_0, \ldots, \sigma_m, 0 \},
\]
and we may define $q: \{ 0, \ldots, m \} \go \{ 0, \ldots, n + 1 \}$ by
$\kappa_{q(i)} = \sigma_i$.  Put
\[
Q(j)		=	q^{-1} \{ 0, \ldots, j \},		
\diagspace
\tau_j		=	(\kappa_j - \kappa_{j + 1}) |Q(j)|,	
\diagspace
\vtr{t}		=	\sum_{j = 0}^n \tau_j \vtr{e}^n_j.
\]
Then $Q \in M(\upr{n}, \upr{m})$, $\vtr{t} \in \Delta^n$, and
$\psi_Q(\vtr{t}) = \vtr{s}$, as required.

For injectivity, let $\vtr{s} \in \Delta^m$ and take $Q': \upr{n'} \gomod
\upr{m}$ and $\vtr{t'} \in \Delta^{n'}$ with $\psi_{Q'} (\vtr{t'}) =
\vtr{s}$.  It will be enough to show that with $\upr{n}$, $Q$ and $\vtr{t}$
as in the previous paragraph, there exists $f: \upr{n} \go \upr{n'}$ such
that
\[
\begin{diagram}
\upr{n}	&			&	\\
\dTo<f	&\rdMod(2,1)^Q		&\upr{m}\\
\upr{n'}&\ruMod(2,1)_{Q'}	&	\\
\end{diagram}
\]
commutes and $f_* \vtr{t} = \vtr{t'}$.  A map $f: \upr{n} \go \upr{n'}$
in $\Dface$ amounts to an $(n + 1)$-element subset of $\upr{n'}$, and it
can be seen that $\{ j' \in \upr{n'} \such t'_j \neq 0 \}$ has $n + 1$
elements.  This defines $f$; the remaining checks are left to the reader.
\end{example*}

\begin{example*}{Edgewise subdivision}	\label{eg:rec-edgewise}
Edgewise subdivision (Figure~\ref{fig:edgewise} and~\cite{Freu})
\begin{figure}
\centering
\setlength{\unitlength}{1mm}
\begin{picture}(36,29)(-18,0)
\cell{0}{3}{b}{\includegraphics{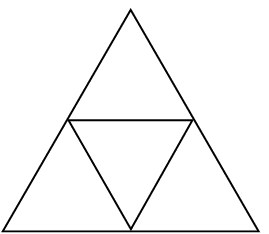}}
\cell{-18}{0}{bl}{00}
\cell{0}{29}{t}{11}
\cell{18}{0}{br}{22}
\cell{-8}{15}{br}{01}
\cell{8}{15}{bl}{12}
\cell{0}{0}{b}{02}
\end{picture}
\hspace*{4em}
\setlength{\unitlength}{1mm}
\begin{picture}(36,29)(-18,0)
\cell{0}{2}{b}{%
\resizebox{24\unitlength}{26\unitlength}{\includegraphics{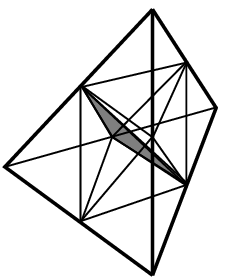}}}
\cell{-13}{12.5}{r}{00}
\cell{4}{2}{t}{11}
\cell{4}{28}{b}{22}
\cell{13}{18}{l}{33}
\cell{9}{10}{l}{13}
\cell{-4}{22}{r}{02}
\cell{-2}{15}{c}{03}
\end{picture}
\caption{Edgewise subdivisions of $\Delta^2$ and $\Delta^3$}
\label{fig:edgewise}
\end{figure}
expresses $\Delta^n$ as $2^n$ smaller copies of itself glued together.
This too can be formalized by a self-similarity structure $M$ on $\Dface$.
Here,
\begin{eqnarray*}
M(\upr{n}, \upr{m})	&
=	&
\{
\textrm{order-preserving injections }
(p, q): \upr{n} \go \upr{m} \times \upr{m}	\\
	&	&
\ \,\textrm{such that }
p(n) \leq q(0)
\}
\end{eqnarray*}
where $\upr{m} \times \upr{m}$ is the product in the category of posets.
For instance, the shaded 2-simplex in Figure~\ref{fig:edgewise} is the
module element $\upr{2} \gomod \upr{3}$ given by the order-preserving
injection $\upr{2} \go \upr{3} \times \upr{3}$ with image $\{ (0, 2), (0,
3), (1, 3) \}$.  Again it can be shown that $\Delta^\blob: \Dface \go
\Top$, equipped with a certain $M$-coalgebra structure, is the universal
solution.
\end{example*}

\section{Spaces admitting a self-similarity}
\label{sec:sass}

Here we answer the question: which topological spaces admit at least one
self-similarity structure?  That is, call a space \demph{self-similar} if it
is homeomorphic to $Ia$ for some self-similarity system $(\scat{A}, M)$
with universal solution $(I, \iota)$ and some $a \in \scat{A}$; which
spaces are self-similar?  We know that every self-similar space is compact
metrizable.  In fact, the converse holds:
\begin{thm}[Self-similar spaces]
\label{thm:selfsimilarspaces}
A topological space is self-similar if and only if it is compact and
metrizable.
\end{thm}
This calls for some explanation.

First recall that when we formalized the notion of self-similarity system
\cite[\S\secselfsimilaritysystems]{SS1}, we allowed there to be infinitely
many equations, even though each individual equation could involve only
finitely many spaces.  So there can be infinite regress: $X_1$ might be
described as a copy of itself copied to a copy of $X_2$, $X_2$ as a copy of
itself glued to a copy of $X_3$, and so on.  This is essential to the proof
of Theorem~\ref{thm:selfsimilarspaces}.  The more restrictive notion of
finite self-similarity is discussed at the end of the section.

Second, this theorem does not exhaust the subject.  It characterizes those
spaces that are self-similar in at least one way, but the same space may
carry several different self-similarities: for instance, the spaces
$\Delta^n$ are self-similar by both barycentric and edgewise
subdivision~(\S\ref{sec:examples}).  Compare the result that every nonempty
set admits at least one group structure, which does not exhaust group
theory.

Here is the strategy for the proof.  Let $S$ be a compact metrizable space.
Cover $S$ by two closed subsets $V_1$ and $V'_1$.  Then $S = V_1 \cup
V'_1$; more precisely, $S$ is the pushout
\[
S = V_1 +_{V''_1} V'_1
\]
where $V''_1 = V_1 \cap V'_1$.  Next, cover $S$ by a different pair $V_2$,
$V'_2$ of closed subsets and write $V''_2 = V_2 \cap V'_2$: then
\begin{eqnarray*}
V_1	&=	&
(V_1 \cap V_2) +_{(V_1 \cap V''_2)} (V_1 \cap V'_2),	\\
V'_1	&=	&
(V'_1 \cap V_2) +_{(V'_1 \cap V''_2)} (V'_1 \cap V'_2),	\\
V''_1	&=	&
(V''_1 \cap V_2) +_{(V''_1 \cap V''_2)} (V''_1 \cap V'_2).	
\end{eqnarray*}
Continue in this way to obtain a countable self-similarity system.  Compact
metrizability of $S$ means that the covers can be chosen to penetrate all
of its structure, and the universal solution $I$ is then made up of the
spaces $S$, $V_1$, $V'_1$, $V''_1$, $V_2$, $V'_2$, $V''_2$, \ldots\ and
inclusions between them.

Formally, a \demph{separating sequence} on a space $S$ is a sequence
$(\mathcal{V}_n)_{n\geq 1}$ of finite closed covers of $S$ such that
\begin{itemize}
\item if $n\geq 1$ and $V, V' \in \mathcal{V}_n$ then $V \cap V' \in
  \mathcal{V}_n$ 
\item if $s, s' \in S$ and for all $n \geq 1$ there exists $V \in
  \mathcal{V}_n$ containing both $s$ and $s'$ then $s = s'$.
\end{itemize}

\begin{lemma}
\label{lemma:separatingsequence}
Every compact metrizable space admits a separating sequence.
\end{lemma}

\begin{proof}
Every compact metrizable space $S$ has a basis $(U_n)_{n\geq 1}$ of open
sets.  For each $n\geq 1$, let
\[
V_n = \textrm{closure}(U_n),
\diagspace
V'_n = S \without U_n,
\diagspace
V''_n = V_n \cap V'_n,
\]
and $\mathcal{V}_n = \{ V_n, V'_n, V''_n \}$.  Then $(\mathcal{V}_n)_{n\geq
1}$ is a sequence of finite closed covers satisfying the binary
intersection axiom.  Let $s$ and $s'$ be distinct points of $S$.  Since $S$
is Hausdorff, there exists $n\geq 1$ with $s \in U_n$ but $s' \not\in V_n$;
then there is no $V \in \mathcal{V}_n$ such that $s, s' \in V$.  \done
\end{proof}

Fix a compact metrizable space $S$ with a separating sequence
$(\mathcal{V}_n)_{n\geq 1}$.

First we define a self-similarity system $(\scat{A}, M)$.  Recall that a
poset can be regarded as a category in which each hom-set has at most one
element: there is a map $a' \go a$ just when $a' \leq a$.  Let $\scat{A} =
\sum_{n\geq 0} \scat{A}_n$, where $\scat{A}_n$ is the sub-poset of
$(\mathcal{V}_1, \sub) \times \cdots \times (\mathcal{V}_n, \sub)$
consisting of those $(V_1, \ldots, V_n)$ for which $V_1 \cap \cdots \cap
V_n \neq \emptyset$.  The module $M$ is also `posetal', that is, $M(b, a)$
never has more than one element; there is a module element
\[
(W_1, \ldots, W_m) 
\gomod
(V_1, \ldots, V_n)
\]
if and only if $m = n+1$ and $W_1 \cap \cdots \cap W_{n+1} \sub V_1 \cap
\cdots \cap V_n$. 

\begin{lemma}
\label{lemma:cmssss}
$(\scat{A}, M)$ is a self-similarity system.
\end{lemma}

\begin{proof}
$M$ is finite since each $\mathcal{V}_n$ is finite.  For nondegeneracy,
the binary intersection axiom guarantees that the poset $\scat{A}$ has all
binary meets.  Hence the category $\scat{A}$ has all pullbacks; they are
the squares in $\scat{A}$ of the form
\[
\begin{diagdiag}
				&	&
(V_1 \cap V'_1, \ldots, V_n \cap V'_n)	&
	&				\\
				&\ldTo	&
					&
\rdTo	&				\\
(V_1, \ldots, V_n)		&	&
					&
	&(V'_1, \ldots, V'_n)		\\
				&\rdTo	&
					&
\ldTo	&				\\
				&	&
(W_1, \ldots, W_n).			&
	&				\\ 
\end{diagdiag}
\]
So by \cite[\corcategorywithpullbacksandequalizers]{SS1}, a functor $X:
\scat{A} \go \Set$ is nondegenerate just when it preserves pullbacks
(existence and preservation of equalizers being trivial).  This is easily
verified when $X = M(b, \dashbk)$ for some $b \in \scat{A}$.  \done
\end{proof}

Next we define an $M$-algebra $(J, \psi)$ in $\Top$.  The functor $J:
\scat{A} \go \Top$ is given by
\[
J(V_1, \ldots, V_n) 
= 
V_1 \cap \cdots \cap V_n
\sub
S, 
\]
and if $V_1 \sub V'_1, \ldots, V_n \sub V'_n$ then $J$ sends the resulting
arrow in $\scat{A}$ to the inclusion $V_1 \cap \cdots \cap V_n \rIncl
V'_1 \cap \cdots \cap V'_n$.  The map $\psi: M \otimes J \go J$ is defined
by taking its component at a module element
\[
(W_1, \ldots, W_{n+1})
\gomod
(V_1, \ldots, V_n),
\]
to be the inclusion
\[
W_1 \cap \cdots \cap W_{n+1} 
\rIncl
V_1 \cap \cdots \cap V_n.
\]

\begin{lemma}	
\label{lemma:fixedpoint}
$J: \scat{A} \go \Top$ is a nondegenerate functor and $\psi: M \otimes J
\go J$ is an isomorphism.
\end{lemma}

\begin{proof}
For the first part, each space $J(V_1, \ldots, V_n)$ is compact Hausdorff,
so it suffices to show that $U \of J: \scat{A} \go \Set$ is nondegenerate;
and as in the proof of~\ref{lemma:cmssss}, this reduces to a
straightforward check that $U \of J$ preserves pullbacks.

We now have to show that for each $(V_1, \ldots, V_n) \in \scat{A}$, the
canonical map
\[
\psi:
\int^{
\begin{scriptarray}
\scriptstyle{(W_1, \ldots, W_{n+1}) \in \scat{A}_{n+1}:}		\\
\scriptstyle{W_1 \cap \cdots \cap W_{n+1} \sub V_1 \cap \cdots \cap V_n}
\end{scriptarray}
}
W_1 \cap \cdots \cap W_{n+1}
\go
V_1 \cap \cdots \cap V_n
\]
is a homeomorphism.  The domain is a finite colimit of compact spaces,
hence compact, and the codomain is Hausdorff, so it suffices to show that
$\psi$ is a bijection.

For surjectivity, let $s \in V_1 \cap \cdots \cap V_n$.  We have $s \in
V_{n+1}$ for some $V_{n+1} \in \mathcal{V}_{n+1}$, so then $s \in V_1 \cap
\cdots \cap V_{n+1}$ with $(V_1, \ldots, V_{n+1}) \in \scat{A}_{n+1}$, and
$\psi$ sends $s \in V_1 \cap \cdots \cap V_{n+1}$ to $s \in V_1 \cap \cdots
\cap V_n$. 

For injectivity, suppose that
\[
t \in W_1 \cap \cdots \cap W_{n+1} \sub V_1 \cap \cdots \cap V_n,
\diagspace
t' \in W'_1 \cap \cdots \cap W'_{n+1} \sub V_1 \cap \cdots \cap V_n
\]
with $\psi(t) = \psi(t')$; then $t$ and $t'$ are equal as elements of $S$.
By the binary intersection property there is a diagram
\begin{equation}	\label{eq:injectivity-diagram}
\begin{diagdiag}
				&		&
(W_1 \cap W'_1, \ldots, W_{n+1} \cap W'_{n+1})	&	
		&				\\
				&\ldTo		&
						&
\rdTo		&				\\
(W_1, \ldots, W_{n+1})		&		&
						&
		&(W'_1, \ldots, W'_{n+1})	\\
\end{diagdiag}
\end{equation}
in $\scat{A}$ and an element 
\[ 
t''
\in 
(W_1 \cap W'_1) \cap \cdots \cap (W_{n+1} \cap W'_{n+1})
\]
equal to $t$ and $t'$ as an element of $S$.
By~\bref{eq:injectivity-diagram}, $t''$ represents the same element of the
coend as both $t$ and $t'$, so $t$ and $t'$ represent the same element of
the coend, as required.  \done
\end{proof}

\begin{propn}
\label{propn:topounivsoln}
$(J, \psi^{-1})$ is the universal solution of $(\scat{A}, M)$ in $\Top$.
\end{propn}

\begin{proof}
Each space $J(V_1, \ldots, V_n)$ is compact, and nonempty by definition of
$\scat{A}$, so it only remains to verify the main part of
condition~\bref{item:pr-notmet} of the Precise Recognition Theorem.  Let $a
= (V_1, \ldots, V_n) \in \scat{A}$, let $t, t' \in V_1 \cap \cdots \cap
V_n$, and suppose there exists $(\cdots \gobymod{m_1} a_0) \in \catI a$
such that $t, t' \in \psi_{m_1} \cdots \psi_{m_k} (Ja_k)$ for all
$k\in\nat$.  Writing $a_k = (V^k_1, \ldots, V^k_{n+k}) \in \scat{A}_{n+k}$,
we have $t, t' \in V^k_1 \cap \cdots \cap V^k_{n+k} \sub V^k_k \in
\mathcal{V}_k$ for all $k\geq 1$, so $t = t'$ by definition of separating
sequence.  \done
\end{proof}

Theorem~\ref{thm:selfsimilarspaces} follows: for if $S$ is nonempty, the
sequence $()$ of length $0$ is an object of $\scat{A}$, and $J() = S$; on
the other hand, the empty space is the universal solution of the
self-similarity system $(\One, \emptyset)$ (Example~\ref{eg:rec-terminal}).

The theorem has a restricted form, as follows.

A self-similarity system $(\scat{A}, M)$ is \demph{discrete} if the
category $\scat{A}$ is discrete (has no arrows except identities).  Any
$\Set$-valued functor on a discrete category is nondegenerate, so a
discrete self-similarity system consists of a set $\scat{A}$ and a family
$(M(b, a))_{b, a \in \scat{A}}$ of sets such that $\sum_b M(b, a)$ is
finite for each $a \in \scat{A}$.  Condition \So\ always holds, so every
discrete self-similarity system has a universal solution: it is $(\oba
\catI, \iota)$, where $\oba\catI : \scat{A} \go \Top$ is the functor
constructed before Lemma~\lemmasequentialcompactness\ of~\cite{SS1}.  To
see this, just note that for each $a \in \scat{A}$, the canonical map
$\pi_a: \ob (\catI a) \go Ia$ is a bijection of compact Hausdorff spaces,
hence a homeomorphism.

A \demph{discretely self-similar} space is one homeomorphic to $Ia$
for some discrete self-similarity system $(\scat{A}, M)$ with universal
solution $(I, \iota)$ and some $a \in \scat{A}$.

\begin{thm}[Discretely self-similar spaces]
\label{thm:discretelyselfsimilarspaces}
The following conditions on a topological space $S$ are equivalent:
\begin{enumerate}
\item	\label{item:dss-dss}
$S$ is discretely self-similar
\item	\label{item:dss-seq} 
$S$ is the limit of some sequence $(\cdots \go S_2 \go S_1)$ of finite
discrete spaces
\item	\label{item:dss-ctbl}
$S$ is the limit of some countable diagram of finite discrete spaces
\item	\label{item:dss-topo}
$S$ is compact, metrizable, and totally disconnected.
\end{enumerate}
\end{thm}

\begin{proof}
\paragraph*{\bref{item:dss-dss}$\implies$\bref{item:dss-seq}}
Let $(\scat{A}, M)$ be a discrete self-similarity system.  The universal
solution is $\oba\catI$, and each $\ob(\catI a)$ is by definition the
limit of a sequence of finite discrete spaces.

\paragraph*{\bref{item:dss-seq}$\implies$\bref{item:dss-ctbl}}
Trivial.

\paragraph*{\bref{item:dss-ctbl}$\implies$\bref{item:dss-topo}}
Compact metrizable spaces are the same as compact Hausdorff second
countable spaces.  The classes of compact Hausdorff spaces and totally
disconnected spaces are closed under all limits, and the class of second
countable spaces is closed under countable limits.

\paragraph*{\bref{item:dss-topo}$\implies$\bref{item:dss-dss}}
$S$ has a basis $(U_n)_{n\geq 1}$ of open sets for which each $U_n$ is also
closed~\cite[Theorem II.4.2]{Joh}.  The separating sequence
$(\mathcal{V}_n)_{n\geq 1}$ constructed in
Lemma~\ref{lemma:separatingsequence} then has the property that each
$\mathcal{V}_n$ is a partition of $S$, so the resulting category $\scat{A}$
is discrete.  This proves the implication when $S$ is nonempty, and the
empty case is trivial.  \done
\end{proof}

\begin{example}
The underlying topological space of the absolute Galois group
$\mr{Gal}(\ovln{\mathbb{Q}}/\mathbb{Q})$ is a countable limit of finite
discrete spaces, so discretely self-similar.
\end{example}

Define \demph{finite self-similarity} as discrete self-similarity
was defined, \latin{mutatis mutandis}.  Since there are
uncountably many homeomorphism classes of compact metrizable
spaces, most self-similar spaces are not finitely self-similar.  

\begin{conj}
The Julia set $J(f)$ of any complex rational function $f$ is finitely
self-similar. 
\end{conj}
This brings us full circle: it says that in the first example of
\cite{SS1}, we could have taken any rational function $f$ and seen the same
type of behaviour: after a finite number of decompositions, no more new
spaces $I_n$ appear.  Both $J(f)$ and its complement are invariant under
$f$, so $f$ restricts to an endomorphism of $J(f)$ and this endomorphism
is, with finitely many exceptions, a $\deg(f)$-to-one mapping.  This
suggests that $f$ itself should provide the self-similarity structure of
$J(f)$, and that if $(\scat{A}, M)$ is the corresponding self-similarity
system then the sizes of $\scat{A}$ and $M$ should be bounded in terms of
$\deg(f)$.

\section{The Cantor set}
\label{sec:thecantorset}

Three closely related classical results (\cite{Wil},~\cite{HY}) describe
the role of the Cantor set among compact metrizable spaces: up to
homeomorphism, 
\begin{itemize}
\item every totally disconnected compact metrizable space is a subspace of
  the Cantor set
\item every nonempty compact metrizable space is a quotient of the Cantor
  set
\item the only nonempty perfect totally disconnected compact metrizable
  space is the Cantor set.
\end{itemize}
These are all consequences of results on self-similarity systems.

We begin by studying discrete self-similarity systems, which, described as
systems of equations \cite[eq.~(\eqinformalsss)]{SS1}, look something like
this:
\begin{eqnarray*}
A	&=	&A + A + B + D	\\
B	&=	&A + C + D	\\
C	&=	&C		\\
D	&=	&D + E		\\
\vdots	&	&\vdots		
\end{eqnarray*}
By adjoining new objects, we can turn the system into one in which each
right-hand side is a sum of at most two terms.  Here, for instance, we can
adjoin objects $A', A'', B', \ldots$ and work with the system
\[
A = A' + D,
\ 
A' = A'' + B, 
\ 
A'' = A + A,
\ 
B = B' + D,
\ 
B' = A + C, 
\ 
\ldots.
\]

\begin{lemma}[$\leq 2$ summands]
\label{lemma:leq2summands}
Let $(\scat{A}, M)$ be a discrete self-similarity system with universal
solution $(I, \iota)$.  Then there is a discrete self-similarity system
$(\scat{A}', M')$ with universal solution $(I', \iota')$ and an injection
$\scat{A} \rIncl \scat{A}'$ such that $\restr{I'}{\scat{A}} \iso I$ and $|
\sum_{b'} M'(b', a') | \leq 2$ for each $a' \in \scat{A}'$.
\end{lemma}

\begin{proof}
Given $a \in \scat{A}$, let $k(a) = | \sum_b M(b, a) |$ and write
\[
\sum_b M(b, a)
=
\{
b^a_1 \gobymod{m^a_1} a,
\ 
\ldots, 
\ 
b^a_{k(a)} \gobymod{m^a_{k(a)}} a
\}.
\]
Let
$
\scat{A}'
=
\{
(a, i)
\such
a \in \scat{A},
\  
0 \leq i \leq k(a)
\}
$
and define $\scat{A} \rIncl \scat{A}'$ by $a \goesto (a, k(a))$.  The
module $M'$ on $\scat{A}'$ has the property that $| M'(b', a') | \leq
1$ for each $b', a' \in \scat{A}'$, and is defined
(Figure~\ref{fig:mprime}) by
\begin{figure}
\begin{diagram}
				&	&
(b^a_1, k(b^a_1))		&	&
(b^a_2, k(b^a_2))		&	&
				&	&
				&	&
(b^a_{k(a)}, k(b^a_{k(a)}))	\\
				&	&
\dMod				&	&
\dMod				&	&
				&	&
				&	&
\dMod				\\
(a, 0)				&\rMod	&
(a, 1)				&\rMod	&
(a, 2)				&\rMod	&
\ 				&\cdots	&
\ 				&\rMod	&
(a, k(a))			\\
\end{diagram}
\caption{Definition of $M'$}
\label{fig:mprime}
\end{figure}
\[
| M'((b, j), (a, i)) | = 1
\iff
i \geq 1
\textrm{ and }
(b, j) 
\in
\{ (a, i - 1), (b^a_i, k(b^a_i)) \}
\]
($b, a \in \scat{A}$, $0 \leq j \leq k(b)$, $0 \leq i \leq k(a)$).  Note
that $(a, i-1) \neq (b^a_i, k(b^a_i))$.  

Let $(I', \iota')$ be the universal solution of $(\scat{A}', M')$ in
$\Top$.  I claim that for all $a \in \scat{A}$ and $i \in \{ 0, \ldots,
k(a) \}$,
\[
I' (a, i)
\iso
I'(b^a_1, k(b^a_1))
+ \cdots +
I'(b^a_i, k(b^a_i)).
\]
This holds when $i = 0$ because
$I' (a, 0) \iso (M' \otimes I') (a, 0) = \emptyset$.  When $i \geq 1$,
\[
I' (a, i)
\iso 
(M \otimes I') (a, i)
=
I' (a, i - 1) + I' (b^a_i, k(b^a_i)),
\]
so the claim follows by induction.  Taking $i = k(a)$, we have in
particular $(\restr{I'}{\scat{A}})a \iso (M \otimes (\restr{I'}{\scat{A}}))
a$ for each $a \in \scat{A}$. 

This argument constructs an isomorphism $\psi: M \otimes
(\restr{I'}{\scat{A}}) \goiso \restr{I'}{\scat{A}}$, described explicitly
as follows.  Any element of the module $M$ is uniquely of the form $b^a_i
\gobymod{m^a_i} a$, and gives rise to a diagram
\begin{equation}	\label{eq:mod-chain}
(b^a_i, k(b^a_i))
\gomod
(a, i)
\gomod 
\cdots
\gomod
(a, k(a))
\end{equation}
in $(\scat{A}', M')$, hence a diagram 
\begin{equation}	\label{eq:map-chain}
I'(b^a_i, k(b^a_i))
\go
I'(a, i)
\go 
\cdots
\go
I'(a, k(a))
\end{equation}
in $\Top$ in which the maps are components of $(\iota')^{-1}$; its
composite is $\psi_{m^a_i}: (\restr{I'}{\scat{A}}) b^a_i \go
(\restr{I'}{\scat{A}}) a$.

We use the Precise Recognition Theorem to show that $(\restr{I'}{\scat{A}},
\psi^{-1})$ is the universal solution of $(\scat{A}, M)$.  It will follow
that $\restr{I'}{\scat{A}} \iso I$.

First, $\restr{I'}{\scat{A}}: \scat{A} \go \Top$ is a nondegenerate functor
sending each $a \in \scat{A}$ to a compact space.  

Second, $\restr{I'}{\scat{A}}$ is occupied.  For take a diagram $\cdots
\gobymod{m_1} a_0$ in $(\scat{A}, M)$.  Then as in~\bref{eq:mod-chain},
each $m_i$ gives rise to a diagram
\[
(a_{i + 1}, k(a_{i + 1}))
\gomod
\cdots
\gomod
(a_i, k(a_i))
\]
in $(\scat{A}', M')$ of length at least $1$.  So there is an infinite
diagram $\cdots \gomod (a_0, k(a_0))$ in $(\scat{A}', M')$.  But $I'$ is
occupied, so $\emptyset \neq I' (a_0, k(a_0)) = (\restr{I'}{\scat{A}})
a_0$.

Third, $| \bigcap_{n \in \nat} \psi_{m_1} \cdots \psi_{m_n} (
(\restr{I'}{\scat{A}}) a_n ) | \leq 1$ for any diagram $\cdots
\gobymod{m_1} a_0$ in $(\scat{A}, M)$.  Indeed, each $\psi_{m_i}$ is the
composite~\bref{eq:map-chain} of one or more components of $(\iota')^{-1}$,
so this follows from the corresponding property of $(\iota')^{-1}$.  \done
\end{proof}

A self-similarity system may contain equations of the form $A = 0$, that
is, there may be objects $a \in \scat{A}$ for which there are no module
elements $\cdot \gomod a$.  If we are only interested in the
self-similarity of nonempty spaces then we can discard those parts of the
system, as follows.

Given a self-similarity system $(\scat{A}, M)$ (not necessarily discrete),
write $\nem{\scat{A}}$ for the full subcategory of $\scat{A}$ whose objects
are those $a \in \scat{A}$ for which $\catI a$ is nonempty, and $\nem{M}$
for the restriction of $M$ to $\nem{\scat{A}}$.

\begin{lemma}[$\geq 1$ summand]
\label{lemma:geq1summand}
Let $(\scat{A}, M)$ be a self-similarity system.  Then $(\nem{\scat{A}},
\nem{M})$ is a self-similarity system.  Moreover, if $(\scat{A}, M)$ has a
universal solution then so does $(\nem{\scat{A}}, \nem{M})$, namely, the
restriction of that of $(\scat{A}, M)$.
\end{lemma}

\begin{proof}
Finiteness of $\nem{M}$ is trivial, and nondegeneracy follows from the
isomorphism $\elt{\nem{M}(b, \dashbk)} \iso \elt{M(b, \dashbk)}$.  So
$(\nem{\scat{A}}, \nem{M})$ is a self-similarity system.

Let $(I, \iota)$ be the universal solution of $(\scat{A}, M)$.  `Moreover'
states that the restriction $(\nem{I}, \nem{\iota})$ of $(I, \iota)$ to
$(\nem{\scat{A}}, \nem{M})$ is the universal solution of $(\nem{\scat{A}},
\nem{M})$.  For this restriction to make sense we need an isomorphism
$\nem{M} \otimes (\restr{I}{\nem{\scat{A}}}) \iso \restr{(M \otimes
I)}{\nem{\scat{A}}}$; indeed, if $a \in \nem{\scat{A}}$ then 
\[
(
\nem{M} \otimes (\restr{I}{\nem{\scat{A}}})
) a
=
\int^{b \in \nem{\scat{A}}}
\nem{M} (b, a) \times Ib
\iso
\int^{b \in \scat{A}}
M(b, a) \times Ib
=
\restr{(M \otimes I)}{\nem{\scat{A}}}
a.
\]
To prove `moreover', just observe that condition~\bref{item:pr-notmet}
(or~\bref{item:pr-met}) of the Precise Recognition Theorem is satisfied by
$(I, \iota)$, hence by $(\nem{I}, \nem{\iota})$.  
\done 
\end{proof}

\begin{cor}[$1$ or $2$ summands]
\label{cor:1or2summands}
Let $S$ be a nonempty discretely self-similar space.  Then there is a
discrete self-similarity system $(\scat{A}, M)$ with universal solution
$(I, \iota)$ and an object $a \in \scat{A}$ such that $S \iso Ia$ and for
each $a' \in \scat{A}$, $\sum_b M(b, a')$ has cardinality $1$ or $2$.
\end{cor}

\begin{proof}
Take a self-similarity system from which $S$ arises and apply
Lemma~\ref{lemma:leq2summands} then Lemma~\ref{lemma:geq1summand}.
\done
\end{proof}

\begin{thm}[Retracts of Cantor set]
\label{thm:retractsofCantorset}
Every nonempty discretely self-similar space is a retract of the Cantor set.
\end{thm}

\begin{proof}
Let $(\scat{A}, M)$ be a discrete self-similarity system such that for each
$a \in \scat{A}$, $\sum_b M(b, a)$ has cardinality $1$ or $2$, and choose
for each $a \in \scat{A}$ a commutative triangle
\[
\begin{diagram}
		&		&2	&		&		\\
		&\ruTo<{s_a}	&	&\rdTo>{r_a}	&		\\
\sum_b M(b, a)	&		&\rTo_1	&		&\sum_b M(b, a) \\
\end{diagram}
\]
in $\Set$.  Write $(I, \iota)$ for the universal solution of $(\scat{A},
M)$.  By the last corollary, it is enough to prove that each $Ia$ is a
retract of the Cantor set.

The Cantor set $2^\nat$, together with the usual isomorphism $\gamma:
2^\nat \goiso 2 \times 2^\nat$, is the universal solution of the
self-similarity system $(\One, 2)$ (Example~\ref{eg:rec-terminal}).  It
gives rise to a functor $\Delta 2^\nat: \scat{A} \go \Top$ constant at
$2^\nat$ and a map
\[
\epsln_a
=
\left(
2^\nat
\rTo^\gamma_\diso
2 \times 2^\nat
\goby{r_a \times 1}
\sum_b M(b, a) \times 2^\nat
=
(M \otimes \Delta 2^\nat) a
\right)
\]
for each $a \in \scat{A}$.  So $(\Delta 2^\nat, \epsln)$ is an
$M$-coalgebra in $\Top$, and there is a unique map $\rho: (\Delta 2^\nat,
\epsln) \go (I, \iota)$ of $M$-coalgebras.

On the other hand, there is a map $\delta: \sum_a Ia \go 2 \times \sum_a
Ia$ whose $a$-component is the composite
\[
Ia 
\rTo^{\iota_a}_\diso
\sum_b M(b, a) \times Ib
\rIncl
\left(
\sum_b M(b, a)
\right) 
\times 
\left(
\sum_b Ib
\right)
\goby{s_a \times 1}
2 \times \sum_b Ib.
\]
There is a unique map $\sigma: (\sum_a Ia, \delta) \go (2^\nat, \gamma)$ of
$2$-coalgebras; preservation of the coalgebra structure says that for each
$a \in \scat{A}$,
\begin{equation}	\label{eq:retract-diag}
\begin{diagram}
Ia		&\rTo^{\delta_a}	&2 \times \sum_b Ib	\\
\dTo<{\sigma_a}	&			&\dTo>{1 \times \sigma}	\\
2^\nat		&\rTo_\gamma		&2 \times 2^\nat	\\
\end{diagram}
\end{equation}
commutes.  It follows that $(\sigma_a)_{a \in \scat{A}}$ is a map $(I,
\iota) \go (\Delta 2^\nat, \epsln)$ of $M$-coalgebras, since for each $a
\in \scat{A}$ there is a commutative diagram as in
Figure~\ref{fig:Cantorcommdiag}.
\begin{figure}
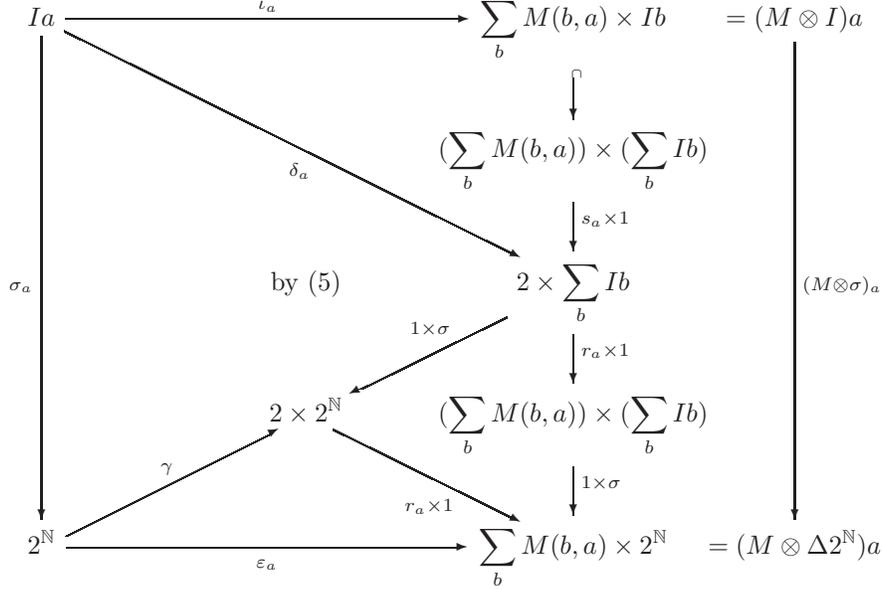

\begin{diagram}[width=5em,height=2.5em]
Ia		&			&\rTo^{\iota_a}		&
			&\sum_b M(b, a) \times Ib		&
= (M \otimes I) a		\\
		&\rdTo(4,4)<{\delta_a}	&			&
			&\dIncl					&\\
		&			&			&
			&(\sum_b M(b, a)) \times (\sum_b Ib)	&\\
		&			&			&
			&\dTo>{s_a \times 1}			&\\
\dTo<{\sigma_a}	&	&\textrm{by \bref{eq:retract-diag}}	&
			&2 \times \sum_b Ib			&
\dTo>{(M \otimes \sigma)_a}	\\
		&			&			&
\ldTo(2,2)<{1 \times \sigma}	&\dTo>{r_a \times 1}		&\\
		&			&2 \times 2^\nat	&
			&(\sum_b M(b, a)) \times (\sum_b Ib)	&\\
		&\ruTo(2,2)<\gamma	&			&
\rdTo(2,2)<{r_a \times 1}&\dTo>{1 \times \sigma}&		&\\
2^\nat		&			&\rTo_{\epsln_a}	&
			&\sum_b M(b, a) \times 2^\nat		&
= (M \otimes \Delta 2^\nat) a	\\
\end{diagram}
\caption{Diagram for the proof of Theorem~\ref{thm:retractsofCantorset}}
\label{fig:Cantorcommdiag}
\end{figure}

We now have maps $(I, \iota) \oppair{\sigma}{\rho} (\Delta 2^\nat,
\epsln)$ of $M$-coalgebras, so by terminality of $(I, \iota)$, the
triangle
\[
\begin{diagram}
	&			&2^\nat	&			&	\\
	&\ruTo<{\sigma_a}	&	&\rdTo>{\rho_a}		&	\\
Ia	&			&\rTo_1	&			&Ia	\\
\end{diagram}
\]
commutes for each $a \in \scat{A}$.
\done
\end{proof}

Using Theorem~\ref{thm:discretelyselfsimilarspaces}, the first classical
result follows:
\begin{cor}
\label{cor:totallydisconnectedretracts}
Every nonempty totally disconnected compact metrizable space is a retract
of the Cantor set.
\done
\end{cor}

We will deduce that every nonempty compact metrizable space is a quotient
of the Cantor set.  This is a topological analogue of the fact that
colimits can be constructed from coproducts and coequalizers, or more
precisely that every finite nonempty colimit is a coequalizer of maps
between finite nonempty coproducts~\cite[V.2]{CWM}.  Recall from
Example~\ref{eg:rec-terminal} that the Cantor set is isomorphic to the
coproduct of $k$ copies of itself for any $k \geq 1$, and is universal as
such for any $k \geq 2$.  Every compact metrizable space is self-similar,
that is, defined by a system of finite colimit diagrams; if we decompose
each colimit as a coequalizer of coproducts, the coproduct part gives the
Cantor set and the coequalizer gives the quotient.

\begin{lemma}
\label{lemma:selfsimilarspaceasquotient}
Every self-similar space is a quotient of a discretely self-similar space.
\end{lemma}
\begin{proof}
Let $(\scat{A}, M)$ be a self-similarity system with universal solution
$(I, \iota)$, let $\scat{A}_0$ be the set of objects of $\scat{A}$
(regarded as a discrete category), and let $M_0$ be the restriction of $M$
to $\scat{A}_0$.  Then the universal solution of $(\scat{A}_0, M_0)$ is
$(\oba\catI, \iota)$, and~\cite[\propnIaasaquotient]{SS1} says that $Ia$ is
a quotient of $(\oba\catI)a$ for each $a \in \scat{A}$.  \done
\end{proof}

\begin{thm}[Quotients of Cantor set]
Every nonempty self-similar space is a quotient of the Cantor set.
\done
\end{thm}

\begin{cor}
Every nonempty compact metrizable space is a quotient of the Cantor set.
\done
\end{cor}

Finally, we deduce the characterization of the Cantor set as the unique
nonempty perfect totally disconnected compact metrizable space.  (A space
is \demph{perfect} if it is not homeomorphic to $T + \{\star\}$ for any
space $T$, or equivalently if it has no isolated points.)

A discrete self-similarity system may contain equations of the form $A =
A$, or loops such as $A = B, B = C, C = A$, or infinite chains such as $A_1
= A_2, A_2 = A_3, \ldots$.  In those cases the universal solution will
involve the one-point space and perhaps other non-perfect spaces too
(Example~\ref{eg:rec-discrete}).  But if the one-point space is not
involved then the universal solution is extremely simple:
\begin{propn}[Empty or Cantor]
\label{propn:emptyorcantor}
Let $(\scat{A}, M)$ be a discrete self-similarity system such that for all
$a \in \scat{A}$, $| Ia | \neq 1$.  Then each $Ia$ is either empty or the
Cantor set.
\end{propn}

\begin{proof}
In Example~\ref{eg:rec-terminal} we saw that for $k\geq 2$, the universal
solution of the self-similarity system $(\One, k)$ is $(2^\nat,
(\psi^{(k)})^{-1} )$.  Here $\psi^{(k)}_p: 2^\nat \go 2^\nat$ is a
contraction with constant $1/3$ for each $p$.  Write $\psi^{(1)}_0$ for the
identity on $2^\nat$.

Define $J: \scat{A} \go \Top$ by
\[
Ja =
\left\{
\begin{array}{ll}
\emptyset	&\textrm{if } Ia = \emptyset	\\
2^\nat		&\textrm{if } Ia \neq \emptyset.
\end{array}
\right.
\]
For each $a \in \scat{A}$, let $k(a) = \left| \sum_{b: Ib \neq \emptyset}
M(b, a) \right|$ and choose an isomorphism between the sets $\sum_{b: Ib
\neq \emptyset} M(b, a)$ and $k(a)$; then for all $a \in \scat{A}$,
\[
(M \otimes J) a
=
\sum_{b \in \scat{A}} M(b, a) \times Jb
\iso
\sum_{b: Ib \neq \emptyset} M(b, a) \times 2^\nat
\iso
k(a) \times 2^\nat.
\]
Define an isomorphism $\gamma_a: Ja \goiso (M \otimes J)a$ for each $a \in
\scat{A}$ as follows.  If $Ia = \emptyset$ then $Ja = \emptyset = (M
\otimes J)a$, so $\gamma_a$ is determined.  If $Ia \neq \emptyset$ then $Ja
= 2^\nat$ and we put
\[
\gamma_a
=
\left(
2^\nat
\rTo^{ ( \psi^{ (k(a)) } )^{-1} }_\diso
k(a) \times 2^\nat
\iso
(M \otimes J) a
\right).
\]

Certainly $(J, \gamma)$ is a fixed point of $M$, $J$ is occupied, and each
$Ja$ is compact and can be equipped with the usual metric.  Now take a
diagram 
\[
\cdots \gobymod{m_2} a_1 \gobymod{m_1} a_0.
\]
For each $r \in \nat$, either $k(a_r) = 1$, in which case
$(\gamma^{-1})_{m_{r+1}}$ is an isometry, or $k(a_r) \geq 2$, in which case
$(\gamma^{-1})_{m_{r+1}}$ is a contraction with constant $1/3$.  The latter
case arises infinitely often: for if not, there is some $s \in \nat$ for
which $1 = k(a_s) = k(a_{s+1}) = \cdots$, and then $|Ia_s| = 1$, contrary
to hypothesis.  Moreover, $\diam(Ja_r) \leq 1$ for each $r$.  So
condition~\bref{item:pr-notmet} of the Precise Recognition Theorem is
satisfied and $(J, \gamma)$ is the universal solution.  \done
\end{proof}

We also need an intuitively obvious lemma: that to calculate the value of
the universal solution at $a$, we can ignore all parts of the system that
do not contribute to $a$.  Precisely, let $(\scat{A}, M)$ be a
self-similarity system and $a \in \scat{A}$.  Write $\scat{A}^{(a)}$ for
the full subcategory of $\scat{A}$ whose objects are those $b \in \scat{A}$
for which there exists a diagram
\begin{equation}	\label{eq:finitetoa}
b = a_n \gobymod{m_n} \cdots \gobymod{m_1} a_0 = a
\end{equation}
for some $n\in\nat$, and $M^{(a)}$ for the restriction of $M$ to
$\scat{A}^{(a)}$.

\begin{lemma}
Let $(\scat{A}, M)$ be a self-similarity system and $a \in \scat{A}$.  Then
$(\scat{A}^{(a)}, M^{(a)})$ is a self-similarity system.  Moreover, if
$(\scat{A}, M)$ has a universal solution then so does $(\scat{A}^{(a)},
M^{(a)})$, namely, the restriction of that of $(\scat{A}, M)$.
\end{lemma}

\begin{proof}
Finiteness of $M^{(a)}$ is trivial.  For nondegeneracy, observe that if $b'
\go b$ in $\scat{A}$ and $b \in \scat{A}^{(a)}$ then $b' \in
\scat{A}^{(a)}$; it follows that $M^{(a)}$ satisfies~\cstyle{ND1}
and~\cstyle{ND2}.  So $(\scat{A}^{(a)}, M^{(a)})$ is a self-similarity
system. 

`Moreover' states that if $(I, \iota)$ is the universal solution of
$(\scat{A}, M)$ then $(I^{(a)}, \iota^{(a)})$ is the universal solution of
$(\scat{A}^{(a)}, M^{(a)})$, where $I^{(a)}$ is the restriction of $I$ and
$\iota^{(a)}$ is the restriction of $\iota$.  For the latter to make sense
we need an isomorphism $\restr{(M \otimes I)}{\scat{A}^{(a)}} \iso M
\otimes \left( \restr{I}{\scat{A}^{(a)}} \right)$; indeed, the square
\[
\begin{diagram}
\ftrcat{\scat{A}}{\Set}		&\rTo^{M \otimes \dashbk}	&
\ftrcat{\scat{A}}{\Set}		\\
\dTo<{\mr{restriction}}		&				&
\dTo>{\mr{restriction}}		\\
\ftrcat{\scat{A}^{(a)}}{\Set}	&\rTo_{M^{(a)} \otimes \dashbk}	&
\ftrcat{\scat{A}^{(a)}}{\Set}	\\
\end{diagram}
\]
commutes up to canonical isomorphism, as if $X \in \ftrcat{\scat{A}}{\Set}$
and $a' \in \scat{A}^{(a)}$ then
\begin{eqnarray*}
\left(
M^{(a)} \otimes \left(
\restr{X}{\scat{A}^{(a)}}
\right)
\right)
a'	&
=	&
\int^{b \in \scat{A}^{(a)}}
M^{(a)} (b, a') \times Xb	\\
	&
\iso	&
\int^{b \in \scat{A}}
M(b, a') \times Xb		\\
	&
=	&
(M \otimes X) a'.
\end{eqnarray*}
To prove `moreover', just observe that condition~\bref{item:pr-notmet}
(or~\bref{item:pr-met}) of the Precise Recognition Theorem is satisfied by
$(I, \iota)$, hence by $(I^{(a)}, \iota^{(a)})$.  
\done
\end{proof}

\begin{thm}[Characterization of Cantor set]
A perfect discretely self-similar space is either empty or the Cantor set.
\end{thm}

\begin{proof}
By the previous lemma we may assume that the space is of the form $Ia$
where $(\scat{A}, M)$ is a discrete self-similarity system with universal
solution $(I, \iota)$, $a \in \scat{A}$, and for all $b \in \scat{A}$ there
exists a finite chain~\bref{eq:finitetoa}.  But
\[
Ia 
\iso
\sum_{a_1} M(a_1, a) \times Ia_1
\iso
\sum_{a_1, a_2} M(a_1, a) \times M(a_2, a_1) \times Ia_2
\iso
\cdots,
\]
so for each $b \in \scat{A}$, $Ib$ is a summand of $Ia$; and since $Ia$ is
perfect, $|Ib| \neq 1$.  The result follows from
Proposition~\ref{propn:emptyorcantor}.  \done
\end{proof}

\begin{cor}
A perfect totally disconnected compact metrizable space is either empty or
the Cantor set.  
\done
\end{cor}

\end{document}